\documentclass[10pt,a4paper]{article}
\usepackage[utf8]{inputenc}
\usepackage[english]{babel}
\usepackage{amsmath}
\usepackage{amsfonts}
\usepackage{amssymb}
\usepackage{mathrsfs}
\usepackage{mathtools}
\usepackage{stmaryrd}
\usepackage{graphicx}
\usepackage{xcolor}
\usepackage{soul}
\usepackage{hyperref}
\DeclareMathOperator{\Spec}{Spec}

\DeclareMathOperator{\cl}{cl}

\DeclareMathOperator{\Gal}{Gal}
\DeclareMathOperator{\length}{length}
\DeclareMathOperator{\disc}{disc}
\DeclareMathOperator{\clo}{cl}

\DeclareMathOperator{\C}{\mathbb{C}}
\DeclareMathOperator{\R}{\mathbb{R}}
\DeclareMathOperator{\Z}{\mathbb{Z}}

\DeclareMathOperator{\Q}{\mathbb{Q}}
\DeclareMathOperator{\F}{\mathbb{F}}
\DeclareMathOperator{\N}{\mathbb{N}}
\DeclareMathOperator{\X}{\mathcal{X}}
\DeclareMathOperator{\Y}{\mathcal{Y}}

\author{Lukas Prader}
\title{An arithmetic zeta function respecting multiplicities}
\date{}

\begin{document}

\maketitle
\footnotetext{The author has been supported by the SFB 1085 ``Higher Invariants''. }

\begin{abstract}
In this paper, we study the arithmetic zeta function
$$\mathscr{Z}_{\X}(s) = \prod_p \prod_{\substack{x \in \X_p \\ \text{closed}}} \Big( \frac{1}{1-|\kappa(x)|^{-s}} \Big)^{\mathfrak{m}_{p}(x)}$$
associated to a scheme $\X$ of finite type over $\Z$, where $\kappa(x)$ denotes the residue field and $\mathfrak{m}_{p}(x)$ the multiplicity of $x$ in $\X_p$. If $\X$ is defined over a finite field, then $\mathscr{Z}_{\X}$ appears naturally in the context of point counting with multiplicities. We prove that $\mathscr{Z}_{\X}$ admits a meromorphic continuation to $\{s \in \C \colon \mathrm{Re}(s) > \dim(\X)-1/2\}$ and determine the order of its pole at $s = \dim(\X)$. Finally, we relate $\mathscr{Z}_{\X}$ to a zeta function $\zeta_f$ encoding the residual factorization patterns of a polynomial $f$.

\end{abstract}

\subsubsection*{Notational conventions} 

We define the set $\N$ of natural numbers by $\N := \{z \in \Z \colon z > 0\}$. \\
For any field $K$, we denote by $\overline{K}$ an algebraic closure of $K$. If $K$ is a number field, then we further denote by $O_K$ its ring of integers. Moreover, if $q$ is a prime power, then $\F_q$ means the (unique) field of cardinality $q$. \\ 
Given a scheme $\X = (\X,\mathcal{O}_{\X})$, we denote by $\mathcal{O}_{\X,x}$ the stalk of $\X$ in $x \in \X$, and by $\mathfrak{m}_{\X,x} \subseteq \mathcal{O}_{\X,x}$ its unique maximal ideal. The residue field $\kappa(x)$ of $x$ is defined by $\kappa(x) := \mathcal{O}_{\X,x}/\mathfrak{m}_{\X,x}$. Further, $\X^{\cl}$ means the subset of closed points in $\X$.\\
For a morphism $\varphi \colon \X \to \Y$ of schemes and a point $y \in \Y$, we denote by  $\X_y := \X \times_{\Y} \Spec(\kappa(y))$ the fibre of $\varphi$ over $y$. Finally, if $\Y = \Spec(R)$ and $R \to S$ is a ring extension, then we shall write $\X \otimes_R S$ instead of $\X \times_{\Spec(R)} \Spec(S)$.

\section{Introduction and motivation}

\subsection{Residual factorization patterns of polynomials}

Throughout this subsection, let $K$ be a number field, and let $n \in \N$. \\
For any polynomial $f \in K[X_1,\ldots,X_n]$, it is possible to find $N \in \N$ such that $f \in O_K[N^{-1}][X_1,\ldots,X_n]$. As a consequence, if $\mathfrak{p} \subseteq O_K[N^{-1}]$ is a maximal ideal (which is the same as a maximal ideal $\mathfrak{p} \subseteq O_K$ avoiding $N$), then we may consider the reduction $\overline{f} \in (O_K/\mathfrak{p})[X_1,\ldots,X_n]$ of $f$ modulo $\mathfrak{p}$.  \\
It is due to Ostrowski \cite{Ost} that $f$ is absolutely irreducible over $K$ (i.e., irreducible over $\overline{K}$) if and only if $\overline{f}$ is absolutely irreducible over $O_K/\mathfrak{p}$ for all but finitely many maximal ideals $\mathfrak{p} \subseteq O_K$. Indeed, this is an easy application of Hilbert's Nullstellensatz, yielding a neat local-global principle. \\
The situation becomes more involved once we try to compare the factorization of $f \in K[X_1,\ldots,X_n]$ with those of $\overline{f} \in (O_K/\mathfrak{p})[X_1,\ldots,X_n]$ for all maximal ideals $\mathfrak{p} \subseteq O_K[N^{-1}]$. As was first observed by Hilbert in \cite{Hil}, even if $f$ is irreducible over $K$, it is not necessarily the case that there exists a single $\mathfrak{p}$ such that $\overline{f}$ is irreducible over $O_K/\mathfrak{p}$. In fact, polynomials with this property are rather special. To make this point, let us restrict for a moment to the case of irreducible and normal\footnote{An irreducible polynomial $f \in \Q[X]$ is called \textit{normal} if the field $\Q[X]/(f)$ is a Galois extension of $\Q$.} polynomials $f \in \Q[X]$, for which we obtain the following characterization: \\
\\
\textbf{Proposition 1.1:} \textit{For an irreducible and normal polynomial $f \in \Q[X]$,  the following assertions are equivalent:} \\

(i) \textit{There exists a prime number $p$ such that $\overline{f} \in \F_p[X]$ is irreducible.} \\

(ii) \textit{$\overline{f} \in \F_p[X]$ is irreducible for infinitely many prime numbers $p$.} \\

(iii) \textit{$\Gal(\mathbb{Q}[X]/(f) \mid \mathbb{Q})$ is cyclic.} \\
\\
\textit{Proof:} Assume that $f \in \Z[N^{-1}][X]$, and let $p \nmid N$ be a prime number. If $p \mid \disc(f)$, then $\overline{f} \in \F_p[X]$ is certainly reducible. Otherwise, \cite[Lem.~4.32, Thm.~4.33]{Nar} suggest that $\overline{f} \in \F_p[X]$ is irreducible if and only if $p$ is inert in $O_K$. However, the latter is possible only if $\Gal(\mathbb{Q}[X]/(f) \mid \mathbb{Q})$ is cyclic, cf. \cite[Cor.~2, p.~367]{Nar}.
Conversely, if $\Gal(\mathbb{Q}[X]/(f) \mid \mathbb{Q})$ is cyclic, then Cebotarev's density theorem guarantees infinitely many inert primes; see \cite[Thm.~7.30]{Nar}. $\hfill \Box$ \\
\\
Nevertheless, it turns out that the number of irreducible factors of $f$ over $K$, counted with multiplicities, is uniquely determined by the factorization patterns of $f$ modulo maximal ideals $\mathfrak{p} \subseteq O_K$. But before giving a precise statement, let us become clear about what we mean by a \textit{factorization pattern}. \\
To this aim, let $\mathfrak{p} \subseteq O_K[N^{-1}]$ be any maximal ideal, and assume that $\overline{f} \in (O_K/\mathfrak{p})[X_1,\ldots,X_n]$ is not a constant polynomial. Then we obtain a unique\footnote{As usually in UFDs, such a factorization is unique up to units and permutation of the factors.} factorization
\begin{equation} \label{eq:fact}
\overline{f} =  u \cdot \overline{f_1}^{\alpha_1} \cdots \overline{f_k}^{\alpha_k}
\end{equation}
for some $k \in \N$, natural numbers $\alpha_1, \ldots, \alpha_k \in \N$, a unit $u \in (O_K/\mathfrak{p})^\times$ and polynomials $f_1,\ldots,f_k \in O_K[X_1,\ldots,X_n]$ with the property that $\overline{f_1}, \ldots, \overline{f_k} \in (O_K/\mathfrak{p})[X_1,\ldots,X_n]$ are irreducible and pairwise relatively prime. Further, for any $1 \leq i \leq k$ we define the \textit{inertia degree} $\delta_i$ of $\overline{f_i}$ to be the number of irreducible factors of $\overline{f_i}$ over $\overline{O_K/\mathfrak{p}}$. 
While the inertia degree of an irreducible univariate polynomial is nothing else than its degree, the former turns out to be the right notion when dealing with multivariate polynomials. Further, the name ``inertia degree'' will be justified by Corollary 2.11. \\
\\
\textbf{Definition 1.2:} \textit{Let $f \in O_K[N^{-1}][X_1,\ldots,X_n]$ and a maximal ideal $\mathfrak{p} \subseteq O_K$ be given. If $N \not\in \mathfrak{p}$ and $\overline{f} \in (O_K/\mathfrak{p})[X_1,\ldots,X_n]$ is not a constant polynomial, then adopting the notation from (\ref{eq:fact}), we define the} factorization pattern \textit{$\mathscr{F}_{\mathfrak{p}}$ of $f$ modulo $\mathfrak{p}$ to be the multiset}
\begin{equation} \label{eq:pat}
\mathscr{F}_{\mathfrak{p}} := \{ (\alpha_1, \delta_1), \ldots, (\alpha_k, \delta_k) \}.
\end{equation}
\textit{In any other case, we set $\mathscr{F}_{\mathfrak{p}} := \emptyset$.} \\
\\
To establish our local-global principle, we mimic the ingenious approach proposed by Riemann in his landmark paper \cite{Rie} and associate to any polynomial $f \in K[X_1,\ldots,X_n]$ a global zeta function $\zeta_f$ constructed as follows: For suitable\footnote{``Suitable'' here expresses the requirement on $\zeta_f(s)$ to converge.} $s \in \C$, we define 
\begin{equation*} \label{eq:eulerprod}
\zeta_f(s) = \prod_{\substack{\mathfrak{p} \subseteq O_K \\ \text{maximal}}} \zeta_{f,\mathfrak{p}}(s)
\end{equation*}
to be an Euler product, where each local factor $\zeta_{f,\mathfrak{p}}$ is supposed to encode the factorization pattern (\ref{eq:pat}) of $f$ modulo $\mathfrak{p}$. 
To this aim, we let
\begin{equation*} \label{eq:eulerfact}
\zeta_{f,\mathfrak{p}}(s) := \prod_{(\alpha,\delta) \in \mathscr{F}_{\mathfrak{p}}} \Big( \frac{1}{1 - \mathrm{N}_{K \mid \Q}(\mathfrak{p})^{- \delta s}}\Big)^{\alpha},
\end{equation*}
where we recall that $\mathrm{N}_{K \mid \Q}(\mathfrak{p}) = |O_K/\mathfrak{p}|$. Note that our definition of $\zeta_f$ depends on the choice of $N$ in Definition 1.2; we shall get rid of this ambiguity in §5 by passing from $\zeta_f$ to a certain equivalence class of zeta functions. \\
Anyway, it is not hard to see that $\zeta_f(s)$ converges for every $s \in \C$ with $\mathrm{Re}(s) > 1$, and that $\zeta_f$ even defines a holomorphic function on this half-plane. This allows us to investigate $\zeta_f$ from an analytic point of view, leading to the following \\
\\
\textbf{Theorem 1.3:} \textit{For every $f \in K[X_1,\ldots,X_n]$, the holomorphic function 
$$\zeta_f \colon \{ s \in \C \colon \mathrm{Re}(s) > 1\} \to \C$$
admits a meromorphic continuation to $\C$. Assume further that $f$ factors as 
$$f = v \cdot F_1^{\beta_1} \cdots F_m^{\beta_m}$$
for some $m \geq 0$, natural numbers $\beta_1, \ldots, \beta_m \in \N$, a unit $v \in K^\times$ and irreducible, pairwise relatively prime polynomials $F_1, \ldots, F_m \in K[X_1,\ldots,X_n]$. Then} 
\begin{equation} \label{eq:T1.3}
\mathrm{ord}_{s=1} \text{ } \zeta_f = - (\beta_1 + \ldots + \beta_m).    
\end{equation}
\\
The whole point is that (\ref{eq:T1.3}) constitutes our desired local-global principle:
Indeed, the right hand side of (\ref{eq:T1.3}) counts (with multiplicities) the irreducible factors of $f$ over $K$, hence is of global nature. On the other hand, the factors $\zeta_{f,\mathfrak{p}}$ and thus $\zeta_f$ itself only depend on the residual factorization patters of $f$, so that the left hand side of (\ref{eq:T1.3}) is of local nature. In this sense, we say that the number of irreducible factors over $K$, counted with multiplicities, is a \textit{local invariant} of $f$. As a special case, we obtain: \\
\\
\textbf{Corollary 1.4:} \textit{A polynomial $f \in K[X_1, \ldots, X_n]$ is irreducible if and only if $\zeta_f$ has a simple pole at $s=1$.}

\subsection{Point counting with multiplicities}

Let $\Y$ be a scheme of finite type over a finite field $\F_q$. Define $N_k := |\Y(\F_{q^k})|$ for every $k \in \N$, where $\Y(\F_{q^k})$ means the set of $\F_q$-morphisms $\Spec(\F_{q^k}) \to \Y$. \\
Note that if $\Y$ is affine, say $\Y = \Spec \big( \F_q[X_1,\ldots,X_n]/(f_1,\ldots,f_m) \big)$ for polynomials $f_1,\ldots,f_m \in \F_q[X_1,\ldots,X_n]$, then $N_k$ is precisely the number of common zeros of $f_1,\ldots,f_m$ in $\mathbb{F}_{q^k}^n$. In general, elements in $\Y(\F_{q^k})$ correspond to pairs $(y,\varphi)$, where $y \in \Y$ is a closed point and $\varphi \colon \kappa(y) \to \mathbb{F}_{q^k}$ an $\F_q$-homomorphism. \\
The \textit{Hasse--Weil zeta function} associated to $\Y$ is defined by 
\begin{equation} \label{eq:HW}
Z(\Y,t) := \exp \Big( \sum_{k = 1}^\infty \frac{N_k}{k} t^k \Big) \in \Q \llbracket t \rrbracket.
\end{equation}
If $\Y$ is a smooth, projective, geometrically connected variety\footnote{By a variety, we mean an integral scheme of finite type over a field.}, then $Z(\Y,t)$ is  subject of the \textit{Weil conjectures}, which were first stated in \cite{Wei}. However, in the current paper, our main interest will be in schemes that are far from being regular. \\
\\
To this aim, recall that the \textit{multiplicity} of a point $y \in \Y$ may be defined as 
\begin{equation} \label{eq:mult}
\mathfrak{m}_{\Y}(y) := \lim_{l \to \infty} d! \frac{\length(\mathcal{O}_{\Y,y}/\mathfrak{m}_{\Y,y}^l)}{l^d} \in \N,
\end{equation}
where $d := \dim(\mathcal{O}_{\Y,y})$. In particular, $y$ is a regular point of $\Y$ if and only if $\mathfrak{m}_{\Y}(y) = 1$. \\
For example, assume that $\Y$ is an affine hypersurface passing through the origin, i.e., $\Y = \Spec \big( \F_q[X_1,\ldots,X_n]/(f) \big)$ for a polynomial $f \in \F_q[X_1,\ldots,X_n] \setminus \{0\}$ with vanishing constant term. Denoting by $y \in \Y$ the closed point corresponding to $0 \in \mathbb{F}_q^n$, then $\mathfrak{m}_{\Y}(y)$ is exactly the minimal degree of a monomial appearing in $f$ with non-zero coefficient. \\
Inspired by (\ref{eq:HW}), we now denote for every $k \in \N$ by 
\begin{equation} \label{eq:Mk}
M_k := \sum_{(y,\varphi) \in \Y(\F_{q^k})} \mathfrak{m}_{\Y}(y)
\end{equation}
the number of pairs $(y,\varphi) \in  \Y(\F_{q^k})$, counted with multiplicities, and consider the zeta function of Hasse--Weil type
$$\widetilde{Z}(\Y,t) := \exp \Big( \sum_{k = 1}^\infty \frac{M_k}{k} t^k \Big) \in \Q \llbracket t \rrbracket.$$
It is a well-known fact that 
\begin{equation} \label{eq:pointcounting}
\zeta_{\Y}(s) := \prod_{y \in \Y^{\clo}} \frac{1}{1-|\kappa(y)|^{-s}} = \exp \Big( \sum_{k = 1}^\infty \frac{N_k}{k} q^{-ks} \Big) = Z(\Y,q^{-s})
\end{equation}
for all $s \in \C$ with $\mathrm{Re}(s) > \dim(\Y)$, which follows by direct computation from 
$$N_k = |\Y(\F_{q^k})| = \sum_{j \mid k} j \cdot |\{ y \in \Y^{\cl} \colon \kappa(y) = \F_{q^j} \}|.$$
Our goal in this subsection is to establish an analogon of (\ref{eq:pointcounting}) for $\widetilde{Z}(\Y,t)$. To do so, we first eliminate any kind of convergence issues by noting that
$M_k \ll_{\Y} N_k$ (i.e., the implied constant depends on $\Y$) as $k \to \infty$, 
which is a consequence of Bezout's theorem (see Lemma 2.6). Hence $\widetilde{Z}(\Y,q^{-s})$ converges for every $s \in \C, \mathrm{Re}(s) > \dim(\Y)$, and letting 
$$a_j := \sum_{\substack{y \in \Y^{\cl}, \\ \kappa(y) = \F_{q^j}}} \mathfrak{m}_{\Y}(y),$$
we conclude that 
$$\widetilde{Z}(\Y,q^{-s}) = \exp \Big( \sum_{k = 1}^\infty \frac{q^{-ks}}{k} \sum_{j \mid k} j a_j\Big) = \exp \Big( \sum_{l=1}^\infty \frac{1}{l} \sum_{j=1}^\infty a_j q^{-jls} \Big)$$
\begin{equation*} \label{eq:pointcountmult}
= \exp \Big( \sum_{j=1}^\infty - a_j \log \big( 1 - q^{-js}\big) \Big) = \prod_{j=1}^\infty \Big( \frac{1}{1-q^{-js}} \Big)^{a_j}
\end{equation*}
$$= \prod_{y \in \Y^{\clo}} \Big( \frac{1}{1-|\kappa(y)|^{-s}} \Big)^{\mathfrak{m}_{\Y}(y)} =: \mathscr{Z}_{\Y}(s).$$
The purpose of this paper is to introduce and study a global version of the zeta function $\mathscr{Z}_{\Y}$ appearing here. Restricting ourselves to the current special case, our elaborations give rise to the following result (see Proposition 4.3): \\
\\
\textbf{Theorem 1.5:} \textit{Let $\Y$ be a scheme of finite type over $\F_q$, and let $D := \dim(\Y)$. Then the Euler product $\mathscr{Z}_{\Y}(s)$ converges for every $s \in \C$ with $\mathrm{Re}(s) > D$. Moreover, $\mathscr{Z}_{\Y}$ admits a meromorphic continuation to the domain 
$$\{s \in \C \colon \mathrm{Re}(s) > D-1/2\}$$
with a pole at $s = D$, whose order is given by the number of $D$-dimensional irreducible components of $\Y$, counted with multiplicities.} \\
\\
By the multiplicity of an irreducible component $\mathcal{Z} \hookrightarrow \Y$ with generic point $\eta \in \mathcal{Z}$, we simply mean the multiplicity $\mathfrak{m}_{\Y}(\eta)$, or equivalently, the length of $\mathcal{O}_{\Y,\eta}$. Moreover, note that $\mathscr{Z}_{\Y}$ actually has an infinite number of poles in $\{s \in \C \colon \mathrm{Re}(s) > D-1/2\}$, all of them located on the line $D + i \R$. A precise description of these poles may be deduced from \cite[Thm.~3]{SerZ}. Finally, we observe: \\
\\
\textbf{Corollary 1.6:} \textit{Let $\Y$ be a scheme of finite type over $\F_q$, and let $D := \dim(\Y)$. Then the number of $D$-dimensional irreducible components of $\Y$, counted with multiplicities, is uniquely determined by the sequence $(M_k)_{k \geq 1}$ introduced in (\ref{eq:Mk}).}

\subsection{Arithmetic zeta functions}

By an \textit{arithmetic scheme}, we mean a scheme of finite type over $\Z$. A global zeta function associated to an arithmetic scheme will thus be called an \textit{arithmetic zeta function}. \\
There is no commonly accepted catalogue of properties that an arithmetic zeta function $\mathfrak{Z}_{\X}$ should satisfy (however, see \cite[Ch.~5]{IK}), but one generally expects it to be given by an Euler product $\mathfrak{Z}_{\X}(s) = \prod_p \mathfrak{Z}_{\X_p}(s)$, where the local factors $\mathfrak{Z}_{\X_p}(s)$ are associated to the fibres $\X_p$ of our given arithmetic scheme $\X$ over prime numbers $p$. In particular, $\mathfrak{Z}_{\X}$ is uniquely determined by the fibres $\X_p$. That means, if a certain invariant of $\X$ is encoded into the analytic behaviour of $\mathfrak{Z}_{\X}$, then this automatically gives rise to a local-global principle. Indeed, this is precisely the way of reasoning that we encountered in §1.1. \\
Moreover, note that the most famous local-global principle of this spirit certainly is the (still open) \textit{Birch--Swinnerton-Dyer conjecture}, where the arithmetic scheme $\X$ is given by an elliptic curve (or more generally, by an abelian variety). \\
\\
In \cite{SerZ}, Serre introduced the arithmetic zeta function 
$$\zeta_{\X}(s) := \prod_{x \in \X^{\cl}} \frac{1}{1-|\kappa(x)|^{-s}} = \prod_p \prod_{x \in \X_p^{\cl}} \frac{1}{1-|\kappa(x)|^{-s}} = \prod_p \zeta_{\X_p}(s),$$
which will thus be referred to as \textit{Serre's zeta function} in the sequel. As we have seen in (\ref{eq:pointcounting}), the local factor $\zeta_{\X_p}$ for a prime $p$ is essentially given by the Hasse--Weil zeta function associated to $\X_p$. Alternatively, one may regard $\zeta_{\X}$ as a product over all closed points in $\X$, where the factor associated to $x \in \X^{\cl}$ only depends on the residue field\footnote{Given a closed point $x \in \X_p \hookrightarrow \X$, note that the residue field of $x$ in $\X_p$ is isomorphic to the residue field of $x$ in $\X$.} $\kappa(x)$ of $x$. \\
We shall summarize the most important analytic properties of $\zeta_{\X}$: \\
\\
\textbf{Theorem 1.7:} \textit{Let $\X$ be an arithmetic scheme of dimension $D := \dim(\X)$. Then the Euler product $\zeta_{\X}(s)$ converges for all $s \in \C$ with $\mathrm{Re}(s) > D$. Moreover, $\zeta_{\X}$ admits a meromorphic continuation to the domain 
$$\{s \in \mathbb{C} \colon \mathrm{Re}(s) > D - 1/2\}$$
having a pole at $s = D$, whose order is given by the number of $D$-dimensional irreducible components of $\X$.} \\
\\
Note further that $s = D$ is the unique pole of $\zeta_{\X}$ precisely if all generic points of $D$-dimensional irreducible components of $\X$ belong to the fibre $\X_0$. Otherwise, there will be infinitely many poles on the line $D + i\R$. We refer to \cite[Thm.~3]{SerZ} for a precise statement. \\ 
\\
In particular, the number of $D$-dimensional irreducible components is a \textit{local invariant} of $\X$ in the sense that it is uniquely determined by the fibres $\X_p$, even by the residue fields of the closed points in $\X_p$ for every prime $p$. \\
Now it is natural to ask whether this local-global principle can be sharpened by considering a refined version of Serre's zeta function $\zeta_{\X}$. \\
Inspired by §1.2, we shall therefore associate to any arithmetic scheme $\X$ the zeta function
$$\mathscr{Z}_{\X}(s) := \prod_p \prod_{x \in \X_p^{\cl}} \Big( \frac{1}{1-|\kappa(x)|^{-s}}\Big)^{\mathfrak{m}_p(x)},$$
where $\mathfrak{m}_{p}(x) := \mathfrak{m}_{\X_p}(x)$ denotes the multiplicity of $x$ in $\X_p$; cf.  (\ref{eq:mult}). Of course, the refinement (in comparison to $\zeta_{\X}$) comes from the ``weights'' $\mathfrak{m}_p(x)$, where the factor in $\mathscr{Z}_{\X}$ corresponding to a closed point $x \in \X_p^{\cl}$ is now allowed to depend on the local ring $\mathcal{O}_{\X_p,x}$ (and not just on the residue field $\kappa(x)$).\\ 
The main result of this paper (stated in a slightly simplified manner) is the forthcoming \\
\\
\textbf{Theorem 1.8:} \textit{Let $\X$ be an arithmetic scheme of dimension $D := \dim(\X)$. Then the Euler product $\mathscr{Z}_{\X}(s)$ converges for all $s \in \C$ with $ \mathrm{Re}(s) > D$. Moreover, $\mathscr{Z}_{\X}$ admits a meromorphic continuation to the domain 
$$\{s \in \C \colon \mathrm{Re}(s) > D-1/2\}$$
with a pole at $s=D$. If the structure morphism $\X \to \Spec(\Z)$ is flat, then this is the unique pole of $\mathscr{Z}_{\X}$, and its order equals the number of $D$-dimensional irreducible components of $\X$, counted with multiplicities.} 
\\
\\
In particular, we obtain the following local-global principle: \\
\\
\textbf{Corollary 1.9:} \textit{Let $\X$ be an arithmetic scheme of dimension $D := \dim(\X)$ such that the structure morphism $\X \to \Spec(\Z)$ is flat. Then the number of $D$-dimensional irreducible components, counted with multiplicities, is a local invariant of $\X$.} \\
\\
Finally, we shall point out a relation between $\mathscr{Z}_{\X^{\mathrm{red}}}$ and Serre's zeta function $\zeta_{\X}$, where $\X^{\mathrm{red}}$ denotes the unique reduced closed subscheme of $\X$ agreeing with $\X$ on the level of points. \\
\\
\textbf{Proposition 1.10:} \textit{Let $\X$ be an arithmetic scheme of dimension $D := \dim(\X)$ such that the structure morphism $\X \to \Spec(\Z)$ is flat. Then the function $s \mapsto \mathscr{Z}_{\X^{\mathrm{red}}}(s)/\zeta_{\X}(s)$ admits a non-vanishing analytic continuation to the domain $\{s \in \C \colon \mathrm{Re}(s) > D-1\}$.} 

\subsection{Structure of this paper}

§2 is the technical heart of the paper and lays the algebro-geometric foundations for our study of $\mathscr{Z}_{\X}$. The latter requires both a qualitative and quantitative understanding of the multiplicities appearing among the closed points of the fibres $\X_p$ for primes $p$. A convenient description of the multiplicities of points in a (single) arithmetic scheme $\X$ is provided by the theory of normal flatness, which was developed by Hironaka in his proof \cite{Hir} for resolution of singularities in characteristic zero. But here we face the additional difficulty that for a point $x \in \X_p$ the quantities $\mathfrak{m}_{\X_p}(x)$ and $\mathfrak{m}_{\X}(x)$ do not necessarily coincide\footnote{For instance, Example 2.2 describes a situation in which $\mathfrak{m}_{\X_p}(x) \geq \mathfrak{m}_{\X}(x)$ holds. On the other hand, given a prime $p$ and a natural number $k$, the irreducible component of $\X = \Spec(\Z[X]/(p^kX))$ corresponding to the principal prime ideal $(p)$ has multiplicity $k$, whereas $\X_p = \Spec(\F_p[X])$ is reduced.}. Nevertheless, we will see that equality holds on an open dense subset of each irreducible component of $\X$ whose generic point belongs to $\X_0$ (cf. Proposition 2.4), and that the sequence $(\mathfrak{m}_{\X_p}(x))_{x,p}$ of multiplicities is bounded (see Proposition 2.7), which suffices for our purpose. For later use, we shall also establish a multivariate version of the famous Kummer--Dedekind theorem \cite[Thm.~4.33]{Nar} (cf. Corollary 2.11).  \\
In §3, we give a brief survey on Dedekind's zeta function $\zeta_K$ associated to a number field $K$ and explain its distinguished role in the study of Serre's zeta function $\zeta_{\X}$. To this aim, we introduce an equivalence relation $\substack{\textcolor{white}{.} \\ \equiv \\ r}$ for meromorphic functions, where $r$ is any non-negative real number, which has a positive impact on the comprehensibility of our statements and proofs. \\ 
Having arrived at §4, we are eventually prepared to investigate $\mathscr{Z}_{\X}$, culminating in Theorem 4.4, which is the main result of this paper and a more precise version of Theorem 1.8. Concerning the proofs, it is worthwhile to consider also a relative version $\mathscr{Z}_{\Y,\X}$ (see Definition 4.1) of our arithmetic zeta function, which may be associated to any closed subscheme $\Y \hookrightarrow \X$. \\
In §5, we shall finally focus on arithmetic schemes of the shape 
$$\X = \Spec \big( O_K[N^{-1}][X_1,\ldots,X_n]/(f) \big)$$
for a polynomial $f \in O_K[N^{-1}][X_1,\ldots,X_n]$, and relate the corresponding zeta function $\mathscr{Z}_{\X}$ to the polynomial zeta function $\zeta_f$ introduced in §1.1. We then conclude the paper by a series of corollaries describing the analytic properties of $\mathscr{Z}_{\X}$ respectively $\zeta_f$ for univariate polynomials $f$.

\subsubsection*{Acknowledgements}

The author is grateful to Herwig Hauser, Guido Kings, Han-Ung Kufner,
Pinaki Mondal and Johannes Sprang for their supportive and fruitful remarks. 

\section{Finiteness results for arithmetic schemes}

Throughout this section, we shall denote by $\X$ an arithmetic scheme such that the structure morphism $\nu \colon \X \to \Spec(\Z)$ is dominant. \\
\\
By a \textit{finiteness result}, we mean the following kind of statement: We are given an invariant $\mathfrak{I}$ for arithmetic schemes, and the claim is that the set of values that $\mathfrak{I}(\X_p)$ may attain as $p$ ranges through all prime numbers is finite. \\
As an initial example, let us consider the following basic \\
\\
\textbf{Lemma 2.1:} \textit{If $\X$ is irreducible, then for every sufficiently large prime number $p$, the fibre $\X_p$ is pure-dimensional of dimension $\dim(\X)-1$.} \\
\\
Concerning the proof, one combines generic flatness with \cite[Cor.~8.2.8]{Liu}. Moreover, the lemma is also true for $p=0$. \\
\\
Note that Lemma 2.1 fits in Grothendieck's philosophy of \textit{spreading out}, which emerges as a special case from the theory developed in \cite[§§8-9]{EGAIV3}: If certain properties are satisfied by $\X_0$, then they are supposed to be immanent in $\X_p$ as well, provided that the prime $p$ is sufficiently large. We will encounter this principle at several points of this section; however, it is not a universal rule: E.g., in the setting of Lemma 2.1, $\X_0$ is irreducible, but it may happen that $\X_p$ fails to be irreducible for every single prime $p$. (Nevertheless, the numbers of irreducible components of $\X_p$ are subject to a plain pattern, to be described in Proposition 2.9.) \\
\\
As already indicated, we desire to gain a good understanding of the multiplicities of (closed) points in the fibres $\X_p$. To this aim, we shall start with a couple of remarks on the concept of multiplicity. \\
Let $(R,\mathfrak{m})$ be a noetherian local ring of dimension $d := \dim(R)$, and define the \textit{Hilbert--Samuel function } $H^{(1)}_R \colon \N \to \N$ of $R$ by $H^{(1)}_R(l) := \length_R (R/\mathfrak{m}^{l+1})$. It is well-known that $H_R^{(1)}(l)$ is, for sufficiently large $l$, given by a polynomial in $l$ of degree $d$ with leading coefficient $e/d!$ for some natural number $e \in \N$, which is defined to be the \textit{multiplicity} of $R$. In particular, if $d = 0$, then $H^{(1)}_R$ is eventually constant and $e = \length_R(R)$. \\
For a locally noetherian scheme $\Y$, the multiplicity $\mathfrak{m}_{\Y}(y)$ of a point $y \in \Y$ is defined to be the multiplicity of $\mathcal{O}_{\Y,y}$, which evidently coincides with the definition given in (\ref{eq:mult}). \\
Moreover, we introduce for any $i \geq 2$ the function $H^{(i)}_R \colon \N \to \N$ given by $H_R^{(i)}(l) = \sum_{j=0}^l H_R^{(i-1)}(j)$. By what we have said, it follows easily that $H_R^{(i)}(l)$ is, for sufficiently large $l$, a polynomial in $l$ of degree $d+i-1$ 
with leading coefficient $e/(d+i-1)!$. 
Consequently, the multiplicity of $R$ may be unravelled from any of the functions $H^{(i)}_R$, where $i \in \N$.\\
\\
Comparing the multiplicities of two Noetherian local rings, the point is that a 
dimensional discrepancy may sometimes be compensated by raising the index of the Hilbert--Samuel function. To give a first taste of this principle, we shall consider the following \\
\\
\textbf{Example 2.2:} Let $\Y$ be a locally noetherian scheme, and let $\mathcal{Z} \hookrightarrow \Y$ be a closed subscheme. Given a point $z \in \mathcal{Z}$, there is a canonical surjection $\mathcal{O}_{\Y,z} \twoheadrightarrow \mathcal{O}_{\mathcal{Z},z}$, and we shall assume that its kernel may be generated by $s \geq 0$ elements.
Then \cite[Thm.~1]{Ben} asserts that \begin{equation} \label{eq:BenT1}
H_{\mathcal{O}_{\Y,z}}^{(i)}(l) \leq H_{\mathcal{O}_{\mathcal{Z},z}}^{(i+s)}(l) 
\end{equation}
for all $i, l \in \N$. If further $\dim(\mathcal{O}_{\mathcal{Z},z}) = \dim(\mathcal{O}_{\Y,z})-s$, then $\deg(H_{\mathcal{O}_{\Y,z}}^{(i)}) = \deg( H_{\mathcal{O}_{\mathcal{Z},z}}^{(i+s)})$ for every $i \in \N$ and hence $\mathfrak{m}_{\Y}(z) \leq \mathfrak{m}_{\mathcal{Z}}(z)$. We shall illustrate this fact by means of two special cases:  \\
(i) Let $p$ be a prime number and $x \in \X_p \hookrightarrow \X$ be any point, then $p$ generates the kernel of $\mathcal{O}_{\X,x} \to \mathcal{O}_{\X_p,x}$. Assuming that all the maximal points of $\X$ belong to $\X_0$, we certainly have $\dim(\mathcal{O}_{\mathcal{X}_p,x}) = \dim(\mathcal{O}_{\X,x})-1$ and thus 
$\mathfrak{m}_{\X_p}(x) \geq \mathfrak{m}_{\X}(x)$. The aim of Proposition 2.4 will be to show that this is ``often'' an equality, in the sense that each irreducible component $\Y \hookrightarrow \X$ has an open dense subscheme $\mathcal{U} \subseteq \Y$ such that $\mathfrak{m}_{\X_p}(x) = \mathfrak{m}_{\X}(x)$ for all $x \in \mathcal{U}_p$ and primes $p$. \\
(ii) Let $K$ be an algebraically closed field, let $f_1,\ldots,f_m \in K[X_1,\ldots,X_n]$ be polynomials with vanishing constant terms, and consider the affine scheme
$$\Y := \Spec \big( K[X_1,\ldots,X_n]/(f_1,\ldots,f_m) \big).$$
Further, denote by $z \in \Y$ the closed point corresponding to $0 \in K^n$. If $\dim(\mathcal{O}_{\Y,z}) > 0$, then it is always possible to find a polynomial $g \in K[X_1,\ldots,X_n]$ of degree one with vanishing constant term so that the closed subscheme 
$$\mathcal{Z} := \Spec \big( K[X_1,\ldots,X_n]/(f_1,\ldots,f_m,g) \big) \hookrightarrow \Y$$
satisfies $\dim(\mathcal{O}_{\mathcal{Z},z}) = \dim(\mathcal{O}_{\Y,z}) -1$. Indeed, for any minimal prime ideal $\mathfrak{p}$ in $\mathcal{O}_{\Y,z}$ we then find $1 \leq j \leq n$ such that $X_j \not\in \mathfrak{p}$. Since $K$ has infinite cardinality, it is now easy to construct a $K$-linear combination $g$ of $X_1,\ldots,X_n$ not contained in any minimal prime ideal of $\mathcal{O}_{\Y,z}$. Moreover, (\ref{eq:BenT1}) suggests that $\mathfrak{m}_{\mathcal{Z}}(z) \geq \mathfrak{m}_{\Y}(z)$. The fact that $g$ may be chosen to be of degree one will become crucial once we are in the context of Bezout's theorem (see Lemma 2.6). $\hfill \diamond$ \\
\\
The functions $H_R^{(i)}$ also play an important role in the theory of \textit{normal flatness}, which we shall briefly discuss now.\\
To this aim, let $\mathcal{Y} \hookrightarrow \mathcal{X}$ be a closed subscheme  defined by a sheaf of ideals $\mathcal{I} \subseteq \mathcal{O}_{\X}$. Following \cite[p.~32]{Ben}, we say that $\X$ is \textit{normally flat} along $\mathcal{Y}$ at a point $y \in \mathcal{Y}$ if $\mathcal{O}_{\Y,y}$ is a regular local ring and $\mathcal{I}_y^l/\mathcal{I}_y^{l+1}$ is a free $\mathcal{O}_{\X,y}/\mathcal{I}_y \cong \mathcal{O}_{\Y,y}$-module for every $l \geq 0$. 
In this case, $\mathcal{I}_{y} \subseteq \mathcal{O}_{\X,y}$ is a prime ideal, and letting $d := \dim(\mathcal{O}_{\Y,y})$, one has 
\begin{equation} \label{eq:normflat}
H^{(i)}_{\mathcal{O}_{\X,y}} = H^{(i+d)}_{(\mathcal{O}_{\X,y})_{\mathcal{I}_{y}}}
\end{equation} 
for every $i \in \N$; a proof can be found in \cite[(2.1.2), p.~33]{Ben}. Since $\X$ is excellent, it follows from \cite[(2.2.1), p.~34]{Ben} that the subset $\mathcal{U} \subseteq \Y$ of all points $y \in \Y$ such that $\X$ is normally flat along $\Y$ at $y$ is open, but possibly empty. However, if $\Y = \overline{\{\eta\}}$ for a point $\eta \in \X$, then $\mathcal{U} \neq \emptyset$ as $\mathcal{O}_{\Y,\eta}$ is a field (see also \cite[(5.2), p.~41]{Ben}), hence it follows from (\ref{eq:normflat}) that 
\begin{equation} \label{eq:multcomp}
\mathfrak{m}_{\X}(y) = \mathfrak{m}_{\X}(\eta) \text{ for every } y \in \mathcal{U}.    
\end{equation}
Proceeding by induction on $\dim(\X)$, one concludes that the set of multiplicities occurring among the points of $\X$ is finite,
and that for any point $x \in \X$ there exists a closed point $z \in \X^{\cl}$ such that $\mathfrak{m}_{\X}(x) = \mathfrak{m}_{\X}(z)$. Of course, the same is true for $\X$ being replaced by any of its fibres $\X_p$. \\
\\
However, the above remarks do not yield an explicit relation between the multiplicities $\mathfrak{m}_{\X}(x)$ and $\mathfrak{m}_{\X_p}(x)$ for any point $x \in \X_p$ and prime $p$. To this aim, the first (qualitative) part of our analysis studies the behaviour of normal flatness under spreading out, culminating in a version of (\ref{eq:multcomp}) that addresses all the fibres $\X_p$ at once. \\
\\
\textbf{Lemma 2.3:} \textit{Let $\mathcal{Y} \hookrightarrow \mathcal{X}$ be a closed subscheme  defined by a sheaf of ideals $\mathcal{I} \subseteq \mathcal{O}_{\X}$ such that the maximal points of $\Y$ belong to the fibre $\Y_0$. Further, denote by $\mathcal{U} \subseteq \mathcal{Y}$ the open subscheme consisting of all points $y \in \Y$ such that $\X$ is normally flat along $\Y$ at $y$. Then $\X_p$ is normally flat along $\Y_p$ at every point in $\mathcal{U}_p$ for all sufficiently large prime numbers $p$.} \\
\\
\textit{Proof:} Since $\X$ is quasi-compact, it suffices to prove the claim in the affine case, say $\X = \Spec(R)$ for a finite type $\Z$-algebra $R$. Then 
$\Y = \Spec(R/I)$ and $\mathcal{I} = \widetilde{I}$ for some ideal $I \subseteq R$ such that $\mathrm{char}(R/I) = 0$. As a consequence, $\X_p = \Spec(\overline{R})$ and $\Y_p = \Spec(\overline{R}/\overline{I})$, where $\overline{R} := R \otimes_{\Z} \mathbb{F}_p \cong R/(p)$ and $\overline{I} := 
I\overline{R}$, so that $\Y_p \hookrightarrow \X_p$ corresponds to the sheaf of ideals $\overline{\mathcal{I}} := \overline{I}^\sim$. For any $l \geq 1$, consider the exact sequence
$$\mathrm{Tor}_1^{\Z}(R/I^l,\F_p) \to I^l \otimes_{\Z} \F_p \to R \otimes_{\Z} \F_p \to (R/I^l) \otimes_{\Z} \F_p \to 0,$$
where $\mathrm{Tor}_1^{\Z}(R/I^l,\F_p) \cong (R/I^l)[p]$. Since $p \not\in I^l$, we have $(R/I^l)[p] \neq 0$ precisely if $p \in R/I^l$ is a zero divisor, or equivalently, contained in a prime ideal associated to $I^l \subseteq R$. But by a result due to Ratliff (see \cite{Rat} or \cite{Bro}), the subset $\bigcup_{l \geq 1} \mathrm{Ass}_R(R/I^l) \subseteq \Spec(R)$ is finite, hence $(R/I^l)[p] =0$ for every $l \geq 1$ provided that $p$ is sufficiently large. At the cost of finitely many primes $p$, we may therefore assume that $I^l \otimes_{\Z} \mathbb{F}_p \cong \overline{I}^l$ for every $l \geq 0$, similarly that $(I^l/I^{l+1}) \otimes_{\Z} \F_p \cong \overline{I}^l/\overline{I}^{l+1}$, and by spreading out that $\mathcal{U}_p \to \Spec(\F_p)$ is smooth. \\
Let now $y \in \mathcal{U}_p$ be arbitrary.  Denoting by $\mathfrak{p} \subseteq R/I$ the prime ideal corresponding to $y \in \Y$, we then know that $\mathcal{I}_y^l/\mathcal{I}_y^{l+1}$ is a free $\mathcal{O}_{\Y,y}$-module for every $l \geq 0$, say 
$$\mathcal{I}_y^l/\mathcal{I}_y^{l+1} \cong \big( I^l/I^{l+1} \big)_{\mathfrak{p}}\cong \bigoplus_{i \in J_l} (R/I)_{\mathfrak{p}}$$
for certain index sets $J_l$. But applying $- \otimes_{\Z} \mathbb{F}_p$, we immediately arrive at 
$$\big( \overline{I}^l/\overline{I}^{l+1} \big)_{\overline{\mathfrak{p}}} \cong \big( I^l/I^{l+1} \big)_{\mathfrak{p}} \otimes_{\Z} \mathbb{F}_p \cong \bigoplus_{i \in J_l} (R/I)_{\mathfrak{p}} \otimes_{\Z} \mathbb{F}_p \cong \bigoplus_{i \in J_l} (\overline{R}/\overline{I})_{\overline{\mathfrak{p}}},$$
proving that $\overline{\mathcal{I}}_y^l / \overline{\mathcal{I}}_y^{l+1} \cong \big( \overline{I}^l/\overline{I}^{l+1} \big)_{\overline{\mathfrak{p}}}$ is a free $\mathcal{O}_{\Y_p,y}$-module for every $l \geq 0$. $\hfill \Box$ \\
\\
\textbf{Proposition 2.4:} \textit{Let $\Y \hookrightarrow \X$ be an irreducible component with generic point $\eta \in \Y$
such that $\nu(\eta) = (0)$. Then there exists a non-empty open subscheme $\mathcal{U} \subseteq \Y$ such that for all $y \in \mathcal{U}_p$ and primes $p$, we have} 
$$\mathfrak{m}_{\X_p}(y) = \mathfrak{m}_{\X}(y) = \mathfrak{m}_{\X}(\eta).$$
\\
Before providing a proof, we need to recall the following lemma, which is a special case of \cite[Prop.~9.8.6]{EGAIV3}. \\
\\
\textbf{Lemma 2.5:} \textit{There is a non-empty open subscheme $\mathcal{W} \subseteq \Spec(\Z)$ such that for all maximal points $\eta \in \X$ with $\nu(\eta) = (0)$ and all $(p) \in \mathcal{W}$, we have 
$$\mathfrak{m}_{\X_p}(\mu) = \length(\mathcal{O}_{\X_p,\mu}) = \length(\mathcal{O}_{\X,\eta}) = \mathfrak{m}_{\X}(\eta)$$
for every maximal point $\mu \in \X_p$ satisfying $\eta \preceq \mu$.} \\
\\
\textit{Proof of Proposition 2.4:} Define $\mathcal{U} \subseteq \Y$ as in Lemma 2.3, and note that $\mathcal{U} \neq \emptyset$ as $\mathcal{O}_{\Y,\eta}$ is a field. Further, let $p$ be a prime number sufficiently large to guarantee that the conclusions of Lemma 2.3 and Lemma 2.5 hold, and that the maximal points in $\Y_p$ are also maximal in $\X_p$. Concerning the latter\footnote{If $\dim(\Y) = \dim(\X)$, which is the only case that will be needed later, then this is immediately clear by Lemma 2.1.}, if $p \not\in \mathfrak{p} \subseteq R$ is a minimal prime ideal, and if $\mathfrak{P} \subseteq R$ is a minimal prime ideal of $\mathfrak{p}+(p)$ but not of $(p)$, then $\mathfrak{P}$ is minimal containing $\mathfrak{p}+\mathfrak{q}$ for some other minimal prime ideal $\mathfrak{q} \subseteq R$. However, there are at most finitely many $\mathfrak{P}$ with this property. \\ 
Let now $y \in \mathcal{U}_p$ be arbitrary, choose $y \in \mathcal{V} \subseteq \X_p$ affine open, say $\mathcal{V} = \Spec(R)$, and write $\Y_p \times_{\X_p} \mathcal{V} = \Spec(R/I)$ for some ideal $I \subseteq R$. 
Denoting by $\mathfrak{q} \in \Spec(R)$ the prime ideal corresponding to $y$, we see that the kernel of $$R_{\mathfrak{q}} \cong \mathscr{O}_{\X_p,y} \twoheadrightarrow \mathscr{O}_{\Y_p,y} \cong (R/I)_{\mathfrak{q}}$$
is precisely $I_{\mathfrak{q}}$, which is hence a prime ideal of $R_{\mathfrak{q}}$ as $y \in \mathcal{U}_p$ is regular. Considering a primary decomposition for $I$, one concludes that $I_{\mathfrak{q}} = \mathfrak{P}_{\mathfrak{q}}$ for some minimal prime ideal $\mathfrak{P} \subseteq R/I$. Denoting by $\mu \in \Y_p$ the maximal point corresponding to $\mathfrak{P}$, we have $(R_{\mathfrak{q}})_{\mathfrak{P}_{\mathfrak{q}}} \cong R_{\mathfrak{P}} \cong \mathscr{O}_{\X_p,\mu}$, so that $\mathfrak{m}_{\X_p}(y) = \mathfrak{m}_{\X_p}(\mu)$ by normal flatness (\ref{eq:normflat}), which in turn equals $\mathfrak{m}_{\X}(\eta)$ by Lemma 2.5. \\
In summary, we have thus shown that the claim of Proposition 2.4 is valid for all but finitely many prime numbers, say for all primes $p > C$, where $C$ is a constant. Replacing $\mathcal{U}$ by its intersection with the open complement of $\bigcup_{p \leq C} \Y_p \subseteq \Y$ completes the proof. $\hfill \Box$  \\
\\
In the second (quantitative) part of our analysis, we shall address the following problem: As mentioned earlier, the set of multiplicities occurring among the points in $\X_p$ is finite for every prime $p$, but a priori, it is not clear whether these multiplicities may grow unboundedly as $p \to \infty$. In fact, by virtue of Bezout's theorem, it is always possible to provide an upper bound that is uniform in $p$. The main ingredient is the following \\
\\
\textbf{Lemma 2.6:} \textit{Let $K$ be a perfect field and let $f_1, \ldots, f_m \in K[X_1,\ldots,X_n] \setminus \{0\}$ be polynomials. Define $C := \max \{\deg(f_1), \ldots, \deg(f_m)\}$. Then any closed point $y$ in 
$$\Y := \Spec \big( K[X_1,\ldots,X_n]/(f_1,\ldots,f_m) \big)$$
has multiplicity $\mathfrak{m}_{\Y}(y) \leq C^n$.} \\
\\
\textit{Proof:} Since $K$ is perfect, we may assume without loss of generality that $K$ is algebraically closed, which is a consequence of (\ref{eq:mult}) and \cite[Prop.~4.7.8]{EGAIV2}. Let $y \in \Y$ be any closed point, and write $d := \dim(\mathcal{O}_{\Y,y})$. Applying Example 2.2 repeatedly, we find degree one polynomials $g_1,\ldots,g_d \in K[X_1,\ldots,X_n]$ so that the closed subscheme
$$\mathcal{Z} := \Spec \big( K[X_1,\ldots,X_n]/(f_1,\ldots,f_m,g_1,\ldots,g_d) \big) \hookrightarrow \Y$$
contains $y$ and satisfies $\dim(\mathcal{O}_{\mathcal{Z},y}) = 0$ as well as $\mathfrak{m}_{\Y}(y) \leq \mathfrak{m}_{\mathcal{Z}}(y)$. In particular, $m+d \geq n$. Assuming that $\deg(f_1) \geq \deg(f_2) \geq \ldots \geq \deg(f_m)$, Bezout's theorem (in the overdetermined case) yields 
$$\mathfrak{m}_{\mathcal{Y}}(y) \leq \mathfrak{m}_{\mathcal{Z}}(y) \leq \prod_{j =1}^{\min\{m,n\}} \deg(f_j) \leq C^n,$$
which is the desired upper bound. $\hfill \Box$ \\
\\
\textbf{Proposition 2.7:} \textit{
There exists a constant $B$ such that $\mathfrak{m}_{\X_p}(x) \leq B$ holds for every point $x \in \X_p$ and prime $p$.} \\
\\
\textit{Proof:} Since $\X$ is quasi-compact, it suffices to consider the affine case. By what we have said, we may further restrict to closed points. We then have
$$\X = \Spec \big( \Z[X_1,\ldots,X_n]/(f_1,\ldots,f_m) \big)$$
for some $f_1,\ldots,f_m \in \Z[X_1,\ldots,X_n]$, hence 
$$\X_p = \Spec \big( \F_p[X_1,\ldots,X_n]/(\overline{f_1},\ldots,\overline{f_m}) \big)$$
for every prime $p$, where $\overline{f_i}$ denotes the reduction of $f_i$ modulo $p$. Now the claim follows from Lemma 2.6. $\hfill \Box$ \\
\\
\textbf{Remark 2.8:} Assume that the structure morphism $\nu \colon \X \to \Spec(\Z)$ factors through $\Spec(O_{K}) \to \Spec(\Z)$ for some number field $K$. In this case, all the finiteness results presented so far are also valid for the fibres of $\X \to \Spec(O_{K})$, which is mainly due to the following simple observation: Given a prime number $p$, and writing $pO_K = \mathfrak{p}_1 \cdots \mathfrak{p}_k$ for maximal ideals $\mathfrak{p}_i \subseteq O_K$, then the fibre $\X_p$ is the disjoint union of the fibres $\X_{\mathfrak{p}_1}, \ldots, \X_{\mathfrak{p}_k}$. In particular, each $\X_{\mathfrak{p}_i}$ is an open subscheme of $\X_p$, so that all our results immediately transfer to this slightly generalized setting. \\
\\
The above elaborations are sufficient for our study of $\mathscr{Z}_{\X}$ in §4. Indeed, Proposition 2.7 plays a key role in the proof that $\mathscr{Z}_{\X}$ is holomorphic, whereas Proposition 2.4 will be applied to establish its meromorphic continuation as well as to describe its pole at $s = \dim(\X)$. \\
For our applications in §5, we shall further discuss the numbers of irreducible components of the fibres $\X_{\mathfrak{p}}$. Even though they do not obey Grothendieck's philosophy of spreading out, it is possible to frame a finiteness result for these numbers; namely, they are determined by the splitting behaviour of the corresponding $\mathfrak{p}$ in a certain number field. Indeed, this fact is crucial for the existence of a meromorphic continuation of Serre's zeta function $\zeta_{\X}$ beyond the line $\dim(\X) +i\R$, as we shall briefly discuss at the very end of §3. \\
Before we start, observe that there is no harm in assuming that $\X$ is reduced, as irreducible components are of purely topological nature. \\
\\
\textbf{Proposition 2.9:} \textit{Suppose that $\X$ is integral with generic point $\eta \in \X$, and that the structure morphism $\nu$ factors through $\Spec(O_K) \to \Spec(\Z)$ for some number field $K$. Moreover, denote by $L$ the algebraic closure\footnote{Note that $L$ is a number field, as $\mathcal{O}_{\X,\eta}$ is a finitely generated field.} of $K$ in $\mathcal{O}_{\X,\eta}$. Then for all but finitely many maximal ideals $\mathfrak{p} \subseteq O_K$, the following holds: If $\mathfrak{p}O_L = \mathfrak{P}_1 \cdots \mathfrak{P}_k$ for maximal ideals $\mathfrak{P}_i \subseteq O_L$ of relative inertia degrees $\delta_i = [O_L/\mathfrak{P}_i : O_K/\mathfrak{p}]$, then the fibre $\X_{\mathfrak{p}}$ has precisely $k$ irreducible components, say $\X_{\mathfrak{p}}^{(1)}, \ldots, \X_{\mathfrak{p}}^{(k)}$, such that $\X_{\mathfrak{p}}^{(i)} \otimes_{O_K/\mathfrak{p}} \overline{O_K/\mathfrak{p}}$ has precisely $\delta_i$ irreducible components for every $1 \leq i \leq k$.}\\
\\
As already indicated, Proposition 2.9 may essentially be extracted from proofs of the meromorphic continuation of $\zeta_{\X}$; see e.g. \cite[Thm.~6.33]{Mus}. Hence we shall allow ourselves to be brief. \\
\\
\textit{Proof:} If the proposition is valid for an integral $O_K$-scheme $\Y$ which is birationally equivalent to $\X$ (i.e., there are open dense $\mathcal{U} \subseteq \X$ and $\mathcal{V} \subseteq \Y$ such that $\mathcal{U} \cong \mathcal{V}$ as $O_K$-schemes), then it must be valid for $\X$ as well: Indeed, by Lemma 2.1, $\mathcal{U}_{\mathfrak{p}}$ (respectively $\mathcal{V}_{\mathfrak{p}}$) contains all the maximal points of $\mathcal{X}_{\mathfrak{p}}$ (respectively $\mathcal{Y}_{\mathfrak{p}}$) for all but finitely many $\mathfrak{p} \subseteq O_K$, and the second claim follows by applying $- \otimes_{O_K/\mathfrak{p}} \overline{O_K/\mathfrak{p}}$ to the irreducible components of $\mathcal{U}_{\mathfrak{p}} \cong \mathcal{V}_{\mathfrak{p}}$. \\
Thus it remains to construct $\Y$ as above for which we are able to verify the claim. As a first step, we may assume that $\X$ is affine. In order to apply resolution of singularities, we now pass from $\X$ to $\X' := \X \otimes_{O_K} K$, but for technical reasons, we further replace $\X'$ by a projective variety $\mathcal{Z}'$ over $K$ containing $\X'$ as an open dense subscheme. Indeed, given a resolution of singularities $\mathcal{Y}' \to \mathcal{Z}'$, this has the advantage that the structure morphism $\mathcal{Y}' \to \Spec(K)$ is not only smooth, but also proper. \\
The latter implies that the global sections $\mathcal{O}_{\mathcal{Y}'}(\mathcal{Y}')$ are a finite field extension of $K$, and since $\mathcal{Y}'$ is regular (hence normal), any germ in $\mathcal{O}_{\mathcal{Y}',\eta} \cong \mathcal{O}_{\mathcal{X},\eta}$ that is algebraic over $K$ must be contained in $\bigcap_{y \in \mathcal{Y}'} \mathcal{O}_{\mathcal{Y}',y} = \mathcal{O}_{\mathcal{Y}'}(\mathcal{Y}')$, hence $\mathcal{O}_{\mathcal{Y}'}(\mathcal{Y}') \cong L$. \\
By spreading out and flat base change, we thus find $N \in \N$ and an integral $O_K$-scheme $\mathcal{Y}$ birationally equivalent to $\X$, equipped with a smooth and proper structure morphism $\Y \to \Spec(O_K[N^{-1}])$, such that $\mathcal{O}_{\Y}(\Y) \cong O_K[N^{-1}]$. Replacing $N$ by a multiple, we may further assume by proper base change  
(see \cite[Thm.~28.1.6]{Vak}) that $\mathcal{O}_{\Y_{\mathfrak{p}}}(\Y_{\mathfrak{p}}) \cong O_L/\mathfrak{p}O_L$ for every maximal ideal $\mathfrak{p} \subseteq O_K[N^{-1}]$. But if $\mathfrak{p}O_L = \mathfrak{P}_1 \cdots \mathfrak{P}_k$ as in the statement of the proposition, then the smoothness and properness of $\Y_{\mathfrak{p}} \to \Spec(O_K/\mathfrak{p})$ suggest that $\Y_{\mathfrak{p}}$ has precisely $k$ irreducible components $\Y_{\mathfrak{p}}^{(1)}, \ldots, \Y_{\mathfrak{p}}^{(k)}$ satisfying $\mathcal{O}_{\Y_{\mathfrak{p}}^{(i)}}(\Y_{\mathfrak{p}}^{(i)}) \cong O_L/\mathfrak{P}_i$. The claim concerning the relative inertia degrees now follows by flat base change $- \otimes_{O_K/\mathfrak{p}} \overline{O_K/\mathfrak{p}}$. $\hfill \Box$ \\
\\
\textbf{Corollary 2.10:} \textit{In the situation of Proposition 2.9, denote by $r_1$ (respectively $r_2$) the number of real (respectively pairs of complex) embeddings of $L$. Then $\X \otimes_{\Z} \R$ respectively $\X \otimes_{\Z} \C$ has precisely $r_1+r_2$ respectively $r_1+2r_2 = [L:\Q]$ irreducible components. In particular, the numbers of irreducible components of $\X \otimes_{\Z} \R$ and $\X \otimes_{\Z} \C$ are uniquely determined by the numbers of irreducible components of $\X_p$ for primes $p$.} \\
\\
\textit{Proof:} We have $L \otimes_{\Q} \R \cong \R^{r_1} \times \C^{r_2}$ and $L \otimes_{\Q} \C \cong \C^{r_1+2r_2} = \C^{[L:\Q]}$. Recalling that $\mathcal{O}_{\Y'}(\Y') \cong L$, the first claim now follows by flat base change $- \otimes_{\Q} \R$ respectively $- \otimes_{\Q} \C$. Concerning the second claim, note that $r_1$ and $r_2$ are uniquely determined by the splitting behaviour of primes in $L$, which is due to \cite[Thm.~1]{Per}. $\hfill \Box$ \\
\\
The Kummer--Dedekind theorem \cite[Thm.~4.33]{Nar} relates the factorization patterns of a univariate irreducible polynomial $f \in \Q[X]$ modulo prime numbers to the splitting behaviour of these primes in the number field $K = \Q[X]/(f)$. However, it seems to be less known that an analogue result is also available for multivariate polynomials: \\
\\
\textbf{Corollary 2.11} (Multivariate Kummer--Dedekind)\textbf{:} \textit{Let $K$ be a number field, and let $f \in O_K[N^{-1}][X_1,\ldots,X_n]$ be a polynomial which is irreducible over $K$. Further, denote by $L$ the algebraic closure of $\Q$ in the field of fractions of $K[X_1,\ldots,X_n]/(f)$. Then for all but finitely many maximal ideals $\mathfrak{p} \subseteq O_K$, the following holds: If $\mathfrak{p}O_L = \mathfrak{P}_1 \cdots \mathfrak{P}_k$ for maximal ideals $\mathfrak{P}_i \subseteq O_L$ of relative inertia degrees $\delta_i = [O_L/\mathfrak{P}_i : O_K/\mathfrak{p}]$, then $\overline{f} \in (O_K/\mathfrak{p})[X_1,\ldots,X_n]$ is a product of precisely $k$ irreducible polynomials in $(O_K/\mathfrak{p})[X_1,\ldots,X_n]$, whose inertia degrees are given by $\delta_1,\ldots,\delta_k$.} \\
\\
Indeed, letting $\X = \Spec \big( O_K[N^{-1}][X_1,\ldots,X_n]/(f) \big)$, then the irreducible components of $\X_{\mathfrak{p}} = \Spec \big( (O_K/\mathfrak{p})[X_1,\ldots,X_n]/(\overline{f}) \big)$ correspond precisely to the irreducible factors of $\overline{f}$ by Krull's principal ideal theorem. Hence Corollary 2.11 is an immediate consequence of Proposition 2.9.

\section{Preliminaries from analysis}

As a starting point, we shall remind the reader of the following well-known \\
\\
\textbf{Lemma 3.1:} \textit{Let $C \geq 0$ be given. Then for any sequence $(a_n)_{n \in \N}$ of real numbers strictly greater than $1$, the following assertions are equivalent:} \\

(i) \textit{$\sum_{n=1}^\infty a_n^{-r}$ converges for every $r > C$.} \\

(ii) \textit{$\phi \colon \{s \in \C \colon \mathrm{Re}(s) > C\} \to \C$ given by 
$$\phi(s) := \prod_{n \in \N} (1-a_n^{-s})$$
defines a non-vanishing holomorphic function. \\
Further, if (i) or (ii) is satisfied, then for every bounded integral sequence $(b_n)_{n \in \N}$, the function $\psi \colon \{s \in \C \colon \mathrm{Re}(s) > C\} \to \C$ defined by 
$$\psi(s) := \prod_{n \in \N} (1-a_n^{-s})^{b_n}$$
is non-vanishing holomorphic as well.} \\
\\
Concerning the proof, it suffices to note that the convergence of $\sum_{n \in \N} a_n^{-r}$ for some $r > 0$ implies the absolute convergence of $\sum_{n \in \N} b_n a_n^{-r}$ by the boundedness of $(b_n)_{n \in \N}$. Note also that the functions $\phi, \psi$ are necessarily non-vanishing on their domains of definition, which is actually built into the definition of convergence for infinite products. \\
In particular, all the zeta functions of our interest are non-vanishing holomorphic on a suitable half-plane. Inspired by this, we shall introduce the following equivalence relation for meromorphic functions, which will make our exposition considerably easier. \\
\\
\textbf{Definition 3.2:} \textit{Let $V \subseteq \C$ be a connected open subset, and let $\phi, \psi \colon V \to \C$ be meromorphic functions which do not vanish identically on $V$. \\
For any real number $r > 0$, we write $\varphi \text{ } \substack{\textcolor{white}{.} \\ \equiv \\ r} \text{ } \psi$ if the function $s \mapsto \phi(s)/\psi(s)$ on $V$ admits a meromorphic continuation to an open subset $W \subseteq \C$ containing  
$U_r := \{ s \in \C \colon \mathrm{Re}(s) > r\}$ such that $\phi/\psi$ is non-vanishing holomorphic on $U_r$. \\
Moreover, we write $\phi \text{ } \substack{\textcolor{white}{.} \\ \equiv \\ 0} \text{ } \psi$ if $\phi/\psi$ is given by a finite product of functions of the shape $(1 - a^{-s})^{\pm 1}$ for $a \in \Z, a \geq 2$.} \\
\\
The point is that if $V$ is an open connected subset of $U_r$ for some $r \geq 0$, and if $\phi, \psi \colon V \to \C$ are meromorphic functions satisfying $\phi \text{ } \substack{\textcolor{white}{.} \\ \equiv \\ r} \text{ } \psi$, then their analytic behaviour must be very similar. For example, if $\phi \colon V \to \C$ is holomorphic, or if $\phi$ has a zero respectively pole at some $s_0 \in V$ of a given order, or if $\phi$ admits a meromorphic continuation to an open subset $V \subseteq V' \subseteq U_r$, then the same must be true for $\psi$. Moreover, if $r = 0$, and if $\phi(z) \in M$ for some integer $z \in \Z \cap V \setminus \{0\}$ and $\Q$-module $M \subseteq \C$, then $\psi(z) \in M$ as well. \\
To put this into context, note that our prototypical example of a non-vanishing holomorphic function on $U_r$ is Serre`s zeta function $\zeta_{\X}$ associated to an arithmetic scheme $\X$ of dimension $\dim(\X) \leq r$ (see Theorem 1.7), and that $\zeta_{\X} \text{ } \substack{\textcolor{white}{.} \\ \equiv \\ 0} \text{ } 1$ if $\X$ is zero-dimensional.\\
For example, let $\X$ be an arithmetic scheme of dimension $D := \dim(\X)$, and let $\mathcal{U} \subseteq \X$ be an open dense subscheme. Then the (reduced) closed subscheme $\mathcal{Y} := \X \setminus \mathcal{U}$ is of dimension strictly less than $D$, thus one has 
$\zeta_{\X} = \zeta_{\mathcal{U}} \cdot \zeta_{\Y} \text{ } \substack{\textcolor{white}{.} \\ \equiv \\ D-1} \text{ } \zeta_{\mathcal{U}}$.    
By what we have said, this means that if one is interested in properties of $\zeta_{\X}$ that are accessible on $U_{D-1}$ (e.g., the pole at $s = D$), then one may replace $\X$ by a suitable open subscheme $\mathcal{U}$, or more generally, by an arithmetic scheme birationally equivalent to $\X$ (in the sense specified in the proof of Proposition 2.9). \\
\\
Now we move on to Dedekind's zeta function, which is a special case of Serre's zeta function, but plays a distinguished role in the investigation of the latter. \\
For a number field $K$, Dedekind's zeta function $\zeta_K$ associated to $K$ is given by 
$$\zeta_K(s) = \zeta_{\Spec(O_K)} =
\prod_{\mathfrak{p}} \frac{1}{1-|O_K/\mathfrak{p}|^{-s}} = \prod_{\mathfrak{p}} \frac{1}{1-\mathrm{N}_{K \mid \Q}(\mathfrak{p})^{-s}}$$
for $s \in \C$ with $\mathrm{Re}(s) > 1$, where the product is over all maximal ideals $\mathfrak{p} \subseteq O_K$. It turns out that $\zeta_K$ has many remarkable properties; we shall gather a few of them in\\
\\
\textbf{Theorem 3.3:} \textit{Let $K$ be a number field. We denote by $r_1$ (respectively $r_2$) the number of real (respectively pairs of complex) embeddings of $K$, and by $\Delta_K$ the discriminant of $K$.} \\

(i) (Meromorphic continuation.) \textit{$\zeta_K$ admits a meromorphic continuation to the whole complex plane $\mathbb{C}$. It is holomorphic on $\mathbb{C} \setminus \{1\}$ with a simple pole at $s = 1$.} \\ 

(ii) (Trivial zeros.) \textit{For every $z \in \N \cup \{0\}$, we have 
\begin{equation*} 
   \mathrm{ord}_{s=-z} \text{ }\zeta_K(s) =
   \begin{cases}
     r_1+r_2-1 & \text{if } z=0, \\
     r_1+r_2 & \text{if } z > 0 \text{ and } 2 \mid z,\\
     r_2 & \text{if } 2 \nmid z.
   \end{cases}
\end{equation*}} 

(iii) (Siegel--Klingen.) \textit{Assume that $K$ is totally real, i.e., that $r_2 = 0$. Then $$\zeta_K(1-z) \in \mathbb{Q} \text{ and } \zeta_K(z) \in \mathbb{Q} \cdot \pi^{z \cdot [K:\mathbb{Q}]} \cdot \Delta_K^{1/2}$$
for every even natural number $z \in 2\mathbb{N}$.} \\

(iv) \textit{Let $L$ be another number field, then $\zeta_K \text{ }  \substack{\textcolor{white}{.} \\ \equiv \\ 0} \text{ } \zeta_L$ implies $\zeta_K = \zeta_L$.} \\
\\
For a proof of (i), we refer to \cite[Thm.~7.3]{Nar}. Then (ii) follows from \cite[Cor.~1, p.~315]{Nar}, the functional equation for $\Phi$ in \cite[Thm.~7.3]{Nar}, and the fact that the Gamma function $\Gamma$ has simple poles at non-positive integers but does not possess any zeros. (iii) is due to \cite{Sie} respectively \cite{Kli}. Finally, (iv) is a consequence of \cite[Thm.~1]{Per}. \\
\\
Now we are prepared to state and discuss a refined version of Theorem 1.7. \\
\\
\textbf{Theorem 3.4:} \textit{Let $\X$ be an integral arithmetic scheme of dimension $D := \dim(\X)$ such that the generic point $\eta \in \X$ belongs to the fibre $\X_0$. Further, denote by $K$ the algebraic closure of $\Q$ in $\mathcal{O}_{\X,\eta}$. Then 
\begin{equation} \label{eq:T3.4}
\zeta_{\X}(s) \text{ }  \substack{\textcolor{white}{.} \\ \equiv \\ D-1/2} \text{ } \zeta_K(s-D+1),
\end{equation}
where $\substack{\textcolor{white}{.} \\ \equiv \\ D-1/2}$ may be replaced by $\substack{\textcolor{white}{.} \\ \equiv \\ 0}$ if $D=1$.} 

\textit{In particular, $\zeta_{\X}$ admits a meromorphic continuation to 
$$\{s \in \C \colon \mathrm{Re}(s) > D - 1/2\}$$
with a unique and simple pole at $s=D$. (If $D=1$, then $\zeta_{\X}$ may even be continued to $\C$, but the pole at $s=D$ is no longer unique in general.)} \\
\\
As a reference, we shall mention \cite[Thm.~9.4]{Ser} and \cite[Thm.~6.33]{Mus}. Note that the appearance of $\zeta_K$ in (\ref{eq:T3.4}) is essentially due to Proposition 2.9, but that the hard part of the proof is to show that $\zeta_{\X}(s)/\zeta_K(s-D+1)$ is non-vanishing holomorphic on $\{s \in \C \colon \mathrm{Re}(s) > D-1/2\}$, which may be deduced e.g. from the Lang--Weil estimates \cite{LW}. \\
Anyway, since $\zeta_K$ admits by Theorem 3.3 a meromorphic continuation to $\C$ with a unique and simple pole at $s=1$, the claim concerning the meromorphic continuation of $\zeta_{\X}$ and its pole is an immediate consequence of (\ref{eq:T3.4}). \\
Given Theorem 3.4, it is easy to obtain an analogue description of $\zeta_{\X}$ for an (almost) arbitrary arithmetic scheme $\X$. We shall not comment on the proof, as a similar argument will be present in Proposition 4.3. \\
\\
\textbf{Corollary 3.5:} \textit{Let $\X$ be an arithmetic scheme of dimension $D := \dim(\X)$ such that the generic points $\eta_1,\ldots,\eta_l \in \X$ of the $D$-dimensional irreducible components $\Y_1,\ldots,\Y_l$ of $\X$ all belong to the fibre $\X_0$. Further, denote for every $1 \leq i \leq l$ by $K_i$ the algebraic closure of $\Q$ in $\mathcal{O}_{\Y,\eta_i}$. Then
$$\zeta_{\X}(s) \text{ }  \substack{\textcolor{white}{.} \\ \equiv \\ D-1} \text{ } \zeta_{\Y_1}(s) \cdots \zeta_{\Y_l}(s) \text{ }  \substack{\textcolor{white}{.} \\ \equiv \\ D-1/2} \text{ } \zeta_{K_1}(s-D+1) \cdots \zeta_{K_l}(s-D+1).$$}

\section{The arithmetic zeta function $\mathscr{Z}_{\X}$}

As a first step in the investigation of $\mathscr{Z}_{\X}$, we shall introduce a relative version of the latter, whose utility will become apparent in the forthcoming proofs. \\
\\
\textbf{Definition 4.1:} \textit{Let $\X$ be an arithmetic scheme, and let $\Y \hookrightarrow \X$ be a closed subscheme. Then we define the arithmetic zeta function $\mathscr{Z}_{\Y,\X}$ of $\Y$ relative to $\X$ by 
$$\mathscr{Z}_{\Y, \X}(s) = \prod_p \prod_{y \in \Y_p^{\cl}} \Big( \frac{1}{1-|\kappa(y)|^{-s}} \Big)^{\mathfrak{m}_{\X_p}(y)}.$$
If $\Y = \X$, we simply write $\mathscr{Z}_{\X}$ instead of $\mathscr{Z}_{\X,\X}$.} \\
\\
Informally speaking, $\mathscr{Z}_{\Y,\X}$ is a product over the closed points in $\Y_p$ for every prime $p$, but with multiplicities coming from $\X_p$. Our main motivation for taking $\mathscr{Z}_{\Y,\X}$ into account is to compensate the fact that $\mathscr{Z}_{\Y}$ and $\mathscr{Z}_{\Y^{\mathrm{red}}}$ may differ drastically (in contrast to Serre's zeta function $\zeta_{\Y} = \zeta_{\Y^{\mathrm{red}}}$). In particular, we even obtain an analogon for $\mathscr{Z}_{\X}$ of the useful identity $\zeta_{\X} = \zeta_{\Y} \cdot \zeta_{\mathcal{U}}$, where $\mathcal{U} \subseteq \X$ denotes the open complement of $\Y$ in $\X$.  \\
\\
\textbf{Lemma 4.2:} \textit{Let $\X$ be an arithmetic scheme, and let $\Y \hookrightarrow \X$ be a closed subscheme of dimension $D := \dim(\Y)$.} \\

(i) \textit{The Euler product $\mathscr{Z}_{\Y,\X}(s)$ converges for every $s \in \C$ with $\mathrm{Re}(s) > D$, so that $\mathscr{Z}_{\Y,\X} \colon \{s \in \C \colon \mathrm{Re}(s) > D\} \to \C$ defines a holomorphic function.} \\

(ii) \textit{Denoting by $\mathcal{U} \subseteq \X$ the open complement of $\Y$ in $\X$, we have $$\mathscr{Z}_{\X} = \mathscr{Z}_{\Y,\X} \cdot \mathscr{Z}_{\mathcal{U}}.$$} 

(iii) \textit{Assume that the structure morphism $\X \to \Spec(\Z)$ factors through $\Spec(O_K) \to \Spec(\Z)$ for some number field $K$. Then 
$$\mathscr{Z}_{\Y, \X}(s) = \prod_{\mathfrak{p} \in \Spec(O_K)^{\cl}} \text{ } \prod_{y \in \mathcal{Y}_{\mathfrak{p}}^{\clo}} \Big( \frac{1}{1-|\kappa(y)|^{-s}} \Big)^{\mathfrak{m}_{\X_\mathfrak{p}}(y)}$$
for every $s \in \C$ with $\mathrm{Re}(s) > D$.} \\
\\
\textit{Proof:} By Proposition 2.7,
the multiplicities $\mathfrak{m}_{\X_p}(x)$ for all $x \in \X_p^{\cl}$ and primes $p$ form a bounded sequence of natural numbers. Hence (i) follows from Theorem 1.7 and Lemma 3.1, (iii) from Remark 2.8. Finally, (ii) is obvious from the definition of $\mathscr{Z}_{\Y,\X}$. $\hfill \Box$ \\
\\
Having established these elementary properties, we shall now decompose $\mathscr{Z}_{\X}$ (with respect to the equivalence relation $\substack{\textcolor{white}{.} \\ \equiv \\ r}$ from Definition 3.2) into a finite product of relative zeta functions $\mathscr{Z}_{\Y,\X}$ associated to certain irreducible components $\Y$ of $\X$, and then study these factors individually. \\
\\
\textbf{Proposition 4.3:} \textit{Let $\X$ be an arithmetic scheme of dimension $D := \dim(\X)$ with structure morphism $\nu \colon \X \to \Spec(\Z)$.} \\

(i) \textit{Denoting by $\Y_1,\ldots,\Y_l$ the $D$-dimensional irreducible components of $\X$, we have}
$$\mathscr{Z}_{\X} \text{ } \substack{\textcolor{white}{.} \\ \equiv \\ D-1} \text{ } \mathscr{Z}_{\mathcal{Y}_1,\X} \cdots \mathscr{Z}_{\mathcal{Y}_l,\X}.$$

(ii) \textit{Let $\Y \hookrightarrow \X$ be a $D$-dimensional irreducible component of $\X$ with generic point $\eta \in \Y$. Then 
$$\mathscr{Z}_{\Y,\X} \text{ } \substack{\textcolor{white}{.} \\ \equiv \\ D-1} \text{ } \zeta_{\Y}^M,$$
where $M = \mathfrak{m}_{\X}(\eta) = \mathfrak{m}_{\X_0}(\eta)$ if $\nu(\eta) = (0)$ and $M = \mathfrak{m}_{\X_p}(\eta)$ if $\nu(\eta) = (p)$ for some prime number $p$.} \\
\\
\textit{Proof:} (i) $\mathscr{Z}_{\X}$ differs from $\mathscr{Z}_{\mathcal{Y}_1,\X} \cdots \mathscr{Z}_{\mathcal{Y}_l,\X}$ by a factor $\xi$ of the shape 
$$\xi(s) = \prod_{z \in \mathcal{Z}^{\cl}} \Big( \frac{1}{1-|\kappa(z)|^{-s}} \Big)^{m_z},$$
where $\mathcal{Z}$ is the (reduced) closed subscheme of $\X$ corresponding to the union of the irreducible components of $\X$ of dimension strictly less than $D$ and of the intersections $\Y_i \cap \Y_j$ of any two distinct $D$-dimensional irreducible components of $\X$, and where $m_z \in \Z$ is an integral multiple of $\mathfrak{m}_{\X_p}(z)$ 
for any $z \in \mathcal{Z}_p^{\cl} \hookrightarrow  \mathcal{Z}^{\cl}$. Note that the sequence $(m_z)_z$ is bounded, which is a consequence of Proposition 2.7 and the fact that the number of irreducible components of $\X$ is finite. Since $\dim(\mathcal{Z}) < D$, the claim thus follows from Theorem 1.7 and Lemma 3.1.  \\
(ii) Assume first that $\nu(\eta) = (p)$. Then $\mathscr{Z}_{\Y,\X} = \mathscr{Z}_{\Y_p,\X_p}$, and as remarked in (\ref{eq:multcomp}), we find an open subscheme $\eta \in \mathcal{V} \subseteq \Y_p$ such that $\mathfrak{m}_{\X_p}(y) = \mathfrak{m}_{\X_p}(\eta)$ for every $y \in \mathcal{V}$. Then $\dim(\Y_p \setminus \mathcal{V}) < D$, hence we arrive at
$$\mathscr{Z}_{\Y_p,\X_p}(s) \text{ } \substack{\textcolor{white}{.} \\ \equiv \\ D-1} \text{ } \prod_{y \in \mathcal{V}^{\clo}} \Big( \frac{1}{1-|\kappa(y)|^{-s}} \Big)^{\mathfrak{m}_{\X_p}(y)} = \zeta_{\mathcal{V}}(s)^{\mathfrak{m}_{\X_p}(\eta)} \text{ } \substack{\textcolor{white}{.} \\ \equiv \\ D-1} \text{ } \zeta_{\Y}(s)^{\mathfrak{m}_{\X_p}(\eta)}.$$
Now suppose that $\nu(\eta) = (0)$. In this case, Proposition 2.4 provides an open subscheme $\eta \in \mathcal{U} \subseteq \Y$ such that $\mathfrak{m}_{\X_p}(y) = \mathfrak{m}_{\X}(\eta)$ for all $y \in \mathcal{U}_p$ and primes $p$. As $\dim(\Y \setminus \mathcal{U}) < D$, we may thus conclude as above that
$$\mathscr{Z}_{\Y,\X}(s) \text{ } \substack{\textcolor{white}{.} \\ \equiv \\ D-1} \text{ } \prod_p \prod_{y \in \mathcal{U}_p^{\clo}} \Big( \frac{1}{1-|\kappa(y)|^{-s}}\Big)^{\mathfrak{m}_{\X_p}(y)} =  \zeta_{\mathcal{U}}(s)^{\mathfrak{m}_{\X}(\eta)} \text{ } \substack{\textcolor{white}{.} \\ \equiv \\ D-1} \text{ } \zeta_{\Y}(s)^{\mathfrak{m}_{\X}(\eta)},$$
completing the proof. $\hfill \Box$ \\
\\
Combining Theorem 3.4, Corollary 3.5 and Proposition 4.3, we obtain  \\
\\
\textbf{Theorem 4.4:} \textit{Let $\X$ be an arithmetic scheme of dimension $D := \dim(\X)$, and suppose that the generic points $\eta_1,\ldots,\eta_l$ of the $D$-dimensional irreducible components $\Y_1,\ldots,\Y_l$ of $\X$ all belong to the fibre $\X_0$. Further, denote for every $1 \leq i \leq l$ by $K_i$ the algebraic closure of $\Q$ in $\mathcal{O}_{\Y_i,\eta_i}$. Then
$$\mathscr{Z}_{\X}(s) \text{ } \substack{\textcolor{white}{.} \\ \equiv \\ D-1} \text{ } 
\zeta_{\Y_1}(s)^{\mathfrak{m}_{\X}(\eta_1)} \cdots \zeta_{\Y_l}(s)^{\mathfrak{m}_{\X}(\eta_l)}
$$
$$\substack{\textcolor{white}{.} \\ \equiv \\ D-1/2} \text{ } \zeta_{K_1}(s-D+1)^{\mathfrak{m}_{\X}(\eta_1)} \cdots \zeta_{K_l}(s-D+1)^{\mathfrak{m}_{\X}(\eta_l)},$$
where $\substack{\textcolor{white}{.} \\ \equiv \\ D-1/2}$ may be replaced by $\substack{\textcolor{white}{.} \\ \equiv \\ 0}$ if $D=1$.}

\textit{In particular, $\mathscr{Z}_{\X}$ admits a meromorphic continuation to 
$$\{ s \in \C \colon \mathrm{Re}(s) > D-1/2 \}$$
with a unique pole at $s = D$, whose order is given by the number of $D$-dimensional irreducible components of $\X$, counted with multiplicities. (If $D=1$, then $\mathscr{Z}_{\X}$ may even be continued to $\C$, but the pole at $s=D$ is no longer unique in general.)} \\
\\
Eventually, we may describe (in a more convenient way) the relation between $\mathscr{Z}_{\X^{\mathrm{red}}}$ and $\zeta_{\X}$ anticipated by Proposition 1.10:  \\
\\
\textbf{Corollary 4.5:} \textit{Let $\X$ be an arithmetic scheme of dimension $D := \dim(\X)$, and suppose that the generic points $\eta_1,\ldots,\eta_l$ of the $D$-dimensional irreducible components $\Y_1,\ldots,\Y_l$ of $\X$ all belong to the fibre $\X_0$. Then
$$\mathscr{Z}_{\X^{\mathrm{red}}} \text{ } \substack{\textcolor{white}{.} \\ \equiv \\ D-1} \text{ } \zeta_{\Y_1} \cdots \zeta_{\Y_l} \text{ } \substack{\textcolor{white}{.} \\ \equiv \\ D-1} \text{ } \zeta_{\X}.$$}

\section{Global zeta classes and local invariants of polynomials}

In this final section, we shall return to the zeta function $\zeta_f$ associated to a polynomial $f \in K[X_1,\ldots,X_n]$ in §1.1, where $K$ is a number field and $n \in \N$.\\
Recall that our definition of $\zeta_f$ depends on the choice of a natural number $N \in \N$ such that $f \in O_K[N^{-1}][X_1,\ldots,X_n]$. However, this problem may be circumvented by replacing $\zeta_f$ by its equivalence class under the equivalence relation $\substack{\textcolor{white}{.} \\ \equiv \\ 0}$ from Definition 3.2.\\
\\
\textbf{Definition 5.1:} \textit{Let $f \in K[X_1,\ldots,X_n]$ be a polynomial, and let $N \in \N$ be such that $f \in O_K[N^{-1}][X_1,\ldots,X_n]$. We define the} global zeta class \textit{$Z_f$ associated to $f$ to be the equivalence class
$$Z_f := [\zeta_f]_{\substack{\equiv \\ 0}}$$
of $\zeta_f$ modulo $\substack{\textcolor{white}{.} \\ \equiv \\ 0}$, which is independent of the choice of $N$.}  \\
\\
Note that $Z_f$ is no longer a well-defined function on $\{s \in \C \colon \mathrm{Re}(s) > 1\}$. 
However, if a certain (analytic) property is valid for any representant of $Z_f$ as soon as it is satisfied by a single representant, then this gives rise to a well-defined property of $Z_f$. Recall that we listed several such properties in the subsequent discussion of Definition 3.2. \\
\\
As a first step, we now follow the idea of Proposition 4.3 and convert the factorization of $f$ over $K$ into a factorization of $Z_f$. \\
\\
\textbf{Proposition 5.2:} \textit{Let $f \in K[X_1,\ldots,X_n]$ be a polynomial.} \\

(i) \textit{Suppose that $f$ factors as 
\begin{equation} \label{eq:P5.2}
f = v \cdot F_1^{\beta_1} \cdots F_m^{\beta_m}
\end{equation}
for some $m \geq 0$, natural numbers $\beta_1, \ldots, \beta_m \in \N$, a unit $v \in K^\times$ and irreducible, pairwise relatively prime polynomials $F_1, \ldots, F_m \in K[X_1,\ldots,X_n]$. Then} 
$$Z_f = Z_{F_1}^{\beta_1} \cdots Z_{F_m}^{\beta_m}.$$

(ii) \textit{If $f$ is irreducible over $K$, then there exists a unique\footnote{Observe that it is not true in general that the number field $L$ is unique with this property. Indeed, there exist non-isomorphic number fields $L, L'$ such that $\zeta_L = \zeta_{L'}$. We refer to \cite{Per} for a detailed discussion of this phenomenon.} Dedekind zeta function $\zeta_L$ associated to a number field $L$ such that $\zeta_L \in Z_f$.} \\
\\
\textit{Proof:} (i) By induction on $m$, it suffices to consider the following two cases: \\
Firstly, assume that $f = F^\beta$ for a polynomial $F \in O_K[N^{-1}][X_1,\ldots,X_n]$ and a natural number $\beta \in \N$. Then it is clear that $\zeta_{f,\mathfrak{p}} = \zeta_{F,\mathfrak{p}}^\beta$ for every maximal ideal $\mathfrak{p} \subseteq O_K[N^{-1}]$, implying that $Z_f = Z_F^\beta$. \\
Secondly, suppose that $f = F \cdot G$ for polynomials $F,G \in O_K[N^{-1}][X_1,\ldots,X_n]$ relatively prime over $K$. Then $\overline{F}, \overline{G} \in (O_K/\mathfrak{p})[X_1,\ldots,X_n]$ will be relatively prime for all but finitely many maximal ideals $\mathfrak{p} \subseteq O_K[N^{-1}]$: Indeed,
assuming to the contrary that $\overline{F}, \overline{G}$ have a common factor over $O_K/\mathfrak{p}$ for infinitely many $\mathfrak{p}$, then it follows from Ostrowski's result \cite{Ost} stated in §1.1 that $F$ and $G$ must have a common factor over $\overline{K}$. However, as $\mathrm{Gal}(\overline{K} \mid K)$ acts transitively on the factors over $\overline{K}$ of any irreducible polynomial over $K$, we see that the common factor of $f$ and $g$ must actually be defined over $K$; a contradiction. \\
Consequently, we have $\zeta_{f,\mathfrak{p}} = \zeta_{F,\mathfrak{p}} \cdot \zeta_{G,\mathfrak{p}}$ for all but finitely many maximal ideals $\mathfrak{p} \subseteq O_K[N^{-1}]$, proving that $Z_f = Z_F \cdot Z_G$. \\
(ii) As in Corollary 2.11, let $L$ be the algebraic closure of $K$ in the field of fractions of $K[X_1,\ldots,X_n]/(f)$. Then for all but finitely many maximal ideals $\mathfrak{p} \subseteq O_K$, if $\mathfrak{p}O_L$ factors as $\mathfrak{p}O_L = \mathfrak{P}_1 \cdots \mathfrak{P}_k$ for maximal ideals $\mathfrak{P}_i \subseteq O_L$ of relative inertia degrees $\delta_i = [O_L/\mathfrak{P}_i : O_K/\mathfrak{p}]$, the local factor $\zeta_{f,\mathfrak{p}}$ is given by 
$$\zeta_{f,\mathfrak{p}}(s) = \prod_{i=1}^k \frac{1}{1-|O_K/\mathfrak{p}|^{- \delta_i s}} = \prod_{i=1}^k \frac{1}{1-p^{- [O_L/\mathfrak{P}_i:\F_p] \cdot s}} = \prod_{i=1}^k \frac{1}{1-\mathrm{N}_{L \mid \Q}(\mathfrak{P}_i)^{-s}},$$
where $p$ is the unique prime number such that $\mathfrak{p} \cap \Z = (p)$. Thus $\zeta_L \in Z_f$ as desired. The uniqueness of $\zeta_L$ follows from Theorem 3.3. $\hfill \Box$ \\
\\
Conversely, if $L$ is any number field, say $L = \Q(\alpha)$ for some $\alpha \in \overline{\Q}$ with minimal polynomial $f \in \Q[X]$ over $\Q$, then the Kummer--Dedekind theorem \cite[Thm.~4.33]{Nar} suggests that $\zeta_L \in Z_f$, proving that any Dedekind zeta function $\zeta_L$ is contained in the global zeta class $Z_f$ of some polynomial $f \in \Q[X]$. In this sense, our global zeta classes may be regarded as a generalization (modulo the relation $\substack{\textcolor{white}{.} \\ \equiv \\ 0}$) of Dedekind's zeta function. \\
Furthermore, note that even though the factorization (\ref{eq:P5.2}) is unique, it is not true in general that $Z_f$ factors \textit{uniquely} into a product $Z_f = Z_{F_1}^{\beta_1} \cdots Z_{F_m}^{\beta_m}$ of global zeta classes $Z_{F_i}$ associated to irreducible polynomials. It seems that this observation (in the context of Dedekind's zeta function) is due to Artin, which led him to the introduction \cite{Art} of what is known today as \textit{Artin's L-function}. \\
\\
Recall that a polynomial invariant is called \textit{local} if it only depends on the residual factorization patterns of the polynomial. \\
\\
\textbf{Corollary 5.3:} \textit{Let $f \in K[X_1,\ldots,X_n]$ be a polynomial. Then $Z_f$ admits a meromorphic continuation to the whole complex plane $\C$. Further, if $f$ factors over $K$ as in (\ref{eq:P5.2}), then 
$$\mathrm{ord}_{s=1} \text{ } Z_f = -(\beta_1 + \ldots + \beta_m).$$
In particular, the number of irreducible factors over $K$, counted with multiplicities, is a local invariant of $f$.} \\
\\
As already indicated, Proposition 5.2 is very much evocative of Proposition 4.3. Hence it does not come as a surprise that there is a close relation between $\zeta_f$ and $\mathscr{Z}_{\X}$ for suitably chosen $\X$. \\
\\
\textbf{Theorem 5.4:} \textit{For any polynomial $f \in K[X_1,\ldots,X_n] \setminus \{0\}$, there exists $N \in \N$ with the following property: We have $f \in O_K[N^{-1}][X_1,\ldots,X_n]$, and defining 
$$\X := \Spec \big( O_K[N^{-1}][X_1,\ldots,X_n]/(f) \big),$$
the relation
$$\mathscr{Z}_{\X}(s) \text{ } \substack{\textcolor{white}{.} \\ \equiv \\ n-1/2} \text{ } \zeta_f(s-n+1)$$
holds. Moreover, if $n=1$, then $\mathscr{Z}_{\X} \in Z_f$.} \\
\\
\textit{Proof:} Assume that $f$ factors over $K$ as in (\ref{eq:P5.2}). Then choose $N \in \N$ in such a way that $v \in O_K[N^{-1}]^\times$ is a unit, the polynomials $F_1,\ldots,F_m$ are defined and irreducible over $O_K[N^{-1}]$, and $O_K[N^{-1}]$ is a principal ideal domain. The latter implies that the maximal points in $\X$ are precisely given by the principal ideals generated by $F_1,\ldots,F_m$. Denoting for every $1 \leq i \leq m$ by $L_i$ the algebraic closure of $\Q$ in the field of fractions of $K[X_1,\ldots,X_n]/(F_i)$, we may conclude from Theorem 4.4 and Proposition 5.2 that 
$$\mathscr{Z}_{\X}(s) \text{ } \substack{\textcolor{white}{.} \\ \equiv \\ n-1/2} \text{ } \zeta_{L_1}(s-n+1)^{\beta_1} \cdots \zeta_{L_m}(s-n+1)^{\beta_m} 
$$
$$\text{ } \substack{\textcolor{white}{.} \\ \equiv \\ n-1} \text{ } \zeta_{F_1}(s-n+1)^{\beta_1} \cdots \zeta_{F_m}(s-n+1)^{\beta_m} \text{ } \substack{\textcolor{white}{.} \\ \equiv \\ n-1} \text{ } \zeta_f(s-n+1).$$
Further, as pointed out in Theorem 4.4, $\substack{\textcolor{white}{.} \\ \equiv \\ n-1/2}$ may be replaced by $ \substack{\textcolor{white}{.} \\ \equiv \\ 0}$ in the above relation if $n=1$. $\hfill \Box$ \\
\\
Finally, we shall present a couple of corollaries describing interesting properties of our global zeta classes $Z_f$, nearby disclosing local invariants of polynomials. For the sake of simplicity, we will state them exclusively for univariate polynomials. (In particular, the corollaries also apply to $\mathscr{Z}_{\X} \in Z_f$ for $\X$ defined as in Theorem 5.4.) However, most of them admit a suitable generalization to the setting of multivariate polynomials. \\
\\
\textbf{Corollary 5.5:} \textit{Assume that $K$ is totally real, and let $f \in K[X] \setminus \{0\}$ be a polynomial. Further, denote by $r$ (respectively $t$) the number\footnote{Note that these numbers are independent of the choice of embedding $K \to \R$.} of real (respectively pairs of conjugate complex) zeros of $f$, counted with multiplicities. Then for every natural number $z \in \N$ we have}
$$\mathrm{ord}_{s=-z} \text{ } Z_f = \begin{cases}
(r+t) \cdot [K:\Q] & \text{if } 2 \mid z,  \\
t \cdot [K:\Q] & \text{if } 2 \nmid z.
\end{cases}$$
\\
\textit{Proof:} Assume that $f$ factors over $K$ as in (\ref{eq:P5.2}), and write $L_i := K[X]/(F_i)$ for every $1 \leq i \leq m$. Denoting by $r_i$ (respectively $t_i$) the number of real (respectively pairs of complex) zeros of $F_i$, then the number field $L_i$ has precisely $r_i \cdot [K:\Q]$ real embeddings and $t_i \cdot [K:\Q]$ pairs of complex embeddings. Thus
$$\mathrm{ord}_{s=-z} \text{ } Z_{F_i} = \mathrm{ord}_{s=-z} \text{ } \zeta_{L_i} = \begin{cases}
(r_i+t_i) \cdot [K:\Q] & \text{if } 2 \mid z,  \\
t_i \cdot [K:\Q] & \text{if } 2 \nmid z
\end{cases}$$
by Theorem 3.3, hence we are done by Proposition 5.2. $\hfill \Box$\\
\\
\textbf{Remark 5.6:} At this point, it is interesting to remark that the splitting field of $f \in K[X]$ is uniquely determined by the set of maximal ideals $\mathfrak{p} \subseteq O_K$ such that $f$ \textit{splits completely} modulo $\mathfrak{p}$, i.e., the factorization pattern $\mathscr{F}_{\mathfrak{p}}$ of $f$ modulo $\mathfrak{p}$ is given by $\deg(f)$ copies of $(1,1)$. Indeed, denoting by $F_1,\ldots,F_m$ the irreducible factors of $f$ over $K$, then up to finitely many exceptional $\mathfrak{p}$, the polynomial $f$ splits completely modulo $\mathfrak{p}$ if and only if $\mathfrak{p}$ splits completely in $K[X]/(F_i)$ for every $1 \leq i \leq m$, which in turn is equivalent by \cite[p.~76]{Mar} to the splitting of $\mathfrak{p}$ in the composite field of the normal closures of the $K[X]/(F_i)$. But the latter is precisely the splitting field of $f$. Now the claim follows from the fact that any finite Galois extension $L$ of a number field $K$ is uniquely determined by the set of maximal ideals in $O_K$ that split completely in $L$; see \cite[Cor.~5, p.~136]{Mar}.  \\
\\
Part (i) of the following corollary was first proved by Schur in \cite{Schur}. \\
\\
\textbf{Corollary 5.7:} \textit{Let $f \in K[X]$ be a polynomial, and denote by $\mathcal{P}_f$ the set of maximal ideals $\mathfrak{p} \subseteq O_K$ such that $f$ has a zero modulo $\mathfrak{p}$.} \\

(i) \textit{$\mathcal{P}_f$ is a finite set if and only if $\deg(f) = 0$ (with the convention that the zero polynomial has degree $- \infty$).} \\

(ii) \textit{If $\deg(f) = d \in \N$ and $\mathrm{ord}_{s=1} \text{ } Z_f < - d/2$, then $f$ has a zero in $K$.} \\
\\
\textit{Proof:} (i) If $\mathcal{P}_f$ is finite, then 
$$\widehat{\zeta_f} := \prod_{\substack{\mathfrak{p} \in \Spec(O_K)^{\cl}, \\ \mathfrak{p} \not\in \mathcal{P}_f}} \zeta_{f,\mathfrak{p}}$$
represents $Z_f$ and has the shape 
$$\widehat{\zeta_f}(s) = \prod_{k = 1}^\infty \frac{1}{1-a_k^{-b_k \cdot s}},$$
where $(a_k)_{k \in \N}$, $(b_k)_{k \in \N}$ are sequences of natural numbers satisfying $b_k \geq 2$ for every $k \in \N$ and  $a_k = \Theta(k \log(k))$ as $k \to \infty$; compare with \cite[Thm.~4.5]{Apo}. Therefore the series
$$\sum_{k=1}^\infty a_k^{-b_k \cdot s}$$
converges absolutely for $s \in \mathbb{C}, \mathrm{Re}(s) > 1/2$, hence the same is true for $\widehat{\zeta_f}(s)$. In particular, this means that $\widehat{\zeta_f}(1) \neq 0$, or equivalently, that $\mathrm{ord}_{s=1} \text{ } Z_f = 0$.\\
(ii) Assuming to the contrary that there is no zero in $K$, and adopting the notation from (\ref{eq:P5.2}), then we have $\deg(F_i) \geq 2$ for every $1 \leq i \leq m$ and therefore $d \geq 2(\beta_1 + \ldots + \beta_m) = -2 \cdot \mathrm{ord}_{s=1} \text{ } Z_f$. $\hfill \Box$ \\
\\
\textbf{Corollary 5.8:} \textit{Let $f \in \mathbb{Q}[X]$ be a polynomial all of whose zeros are real. Assume further that 
$$f = v \cdot F_1^{\beta_1} \cdots F_m^{\beta_m}$$
as in (\ref{eq:P5.2}). Then for every even natural number $z \in 2\mathbb{N}$, we have
$$Z_f(1-z) \in \mathbb{Q} \text{ and } Z_f(z) \in \mathbb{Q} \cdot \pi^{z \cdot \deg(f)} \cdot \disc \Big( \prod_{\substack{ 1 \leq i \leq m, \\ 2 \nmid \beta_i}} F_i \Big)^{1/2}.$$} 
\\
\textit{Proof:} Writing $K_j := \mathbb{Q}[X]/(F_i)$ for every $1 \leq j \leq m$, then $$Z_{F_i}(1-z) \in \mathbb{Q} \text{ and }
Z_{F_i}(z) \in \mathbb{Q} \cdot \pi^{z \cdot [K_j : \mathbb{Q}]} \Delta_{K_j}^{1/2} = \mathbb{Q} \cdot \pi^{z \cdot \deg(F_j)} \disc(F_j)^{1/2}$$
by Theorem 3.3 and the fact that $\disc(F_i)/\Delta_{K_i} \in (\Q^\times)^2$. Observing that 
$$\disc(G \cdot H) = \disc(G) \cdot \disc(H) \cdot \mathrm{Res}(G,H)^2$$ for any two non-constant polynomials $G,H \in \Q[X]$, the claim now follows from Proposition 5.2. $\hfill \Box$ \\
\\
To conclude this paper, we shall demonstrate by means of an example that the discriminant of a polynomial is not a local invariant. \\
\\
\textbf{Example 5.9:} Let $f,g \in \mathbb{Z}[X]$ be given by 
$$f = X^3+6X^2+9X+1 \text{ and } g = X^3+18X^2+81X+27.$$ 
Then $f$ and $g$ are irreducible and $\mathrm{disc}(f)= 3^4 \neq 3^{10} = \mathrm{disc}(g)$. Nevertheless, the factorization patterns of $f$ and $g$ coincide modulo every prime number. To see this, observe that any zero $\alpha \in \overline{\Q}$ of $f$ satisfies $g(3\alpha) = 0$, proving that
$\mathbb{Q}[X]/(f) \cong \mathbb{Q}[X]/(g)$. By \cite[Lem.~4.32, Thm.~4.33]{Nar}, this already implies that the factorization patterns of $f$ and $g$ agree modulo every prime number $p \neq 3$. Thus we are done by noting that
$$\overline{f} = X^3 + \overline{1} = (X+\overline{1})^3 \in \F_3[X] \text{ and } \overline{g} = X^3 \in \F_3[X].$$

\textit{Postal address:} Lukas Prader, Fakultät für Mathematik, Universität Regensburg, 93040 Regensburg, Germany \\

\textit{E-mail address:} lukas.prader@ur.de

\end{document}